\documentclass[11pt]{amsart}

\usepackage{  
  amsmath,
  amsfonts,
  amssymb,
  mvsymb,   
  epsf, 
  epic, 
  eepic, 
  graphicx, 
  picins,
  times,
  color, 
  hyphenat}

\input xy
\xyoption{all}

\newcommand{\brak}[1]{\langle #1\rangle}

\newtheorem{thm}{Theorem}[section]
\newtheorem{lem}[thm]{Lemma}
\newtheorem{cor}[thm]{Corollary}
\newtheorem{prop}[thm]{Proposition}

\newtheorem{rem}{Remark}
\theoremstyle{definition}
\newtheorem{defn}[thm]{Definition}

\newcommand{\bZ}{\mathbb{Z}}
\newcommand{\bQ}{\mathbb{Q}}

\newcommand{\bC}{\mathbb{C}}

\newcommand{\cU}{\mathcal{U}}

\DeclareMathOperator{\Hom}{Hom}

\newcommand{\ra}{\rightarrow}

\newcommand{\kh}[3]{\mbox{KH}^{{#1},{#2}}({#3},\bC)}

\newcommand{\lk}{\mbox{{\em lk}}}

\newcommand{\uhtgen}[2]{U_{a,b,c}^{{#1}}({#2})}
\newcommand{\uht}[2]{U_{a,b,c}^{{#1}}({#2},\bC)}

\newcommand{\uhtabc}[5]{U_{{#1},{#2},{#3}}^{{#4}}({#5},\bC)}
\newcommand{\cut}[2]{C^{{#1}}_{a,b,c}({#2},\bC)}

\title{The universal $sl_3$-link homology}
\author{Marco Mackaay}
\address{Departamento de Matem\'{a}tica\\ Universidade do Algarve\\ 
Campus de Gambelas\\ 8005-139 Faro\\ Portugal}
\email{mmackaay@ualg.pt}
\author{Pedro Vaz}
\address{Departamento de Matem\'{a}tica\\ Universidade do Algarve\\ 
Campus de Gambelas\\ 8005-139 Faro\\ Portugal}
\email{pfortevaz@ualg.pt}

\begin{document}


\begin{abstract}
We define the universal $sl_3$-link homology, which depends on 3 parameters, 
following Khovanov's approach with foams. We show that this 3-parameter 
link homology, when taken with complex coefficients, can be divided 
into 3 isomorphism classes. The first class is 
the one to which Khovanov's original $sl_3$-link homology belongs, the second 
is the one studied by Gornik in the context of matrix factorizations 
and the last one is new. Following an approach similar to Gornik's we show 
that this new link homology can be described in terms of Khovanov's original 
$sl_2$-link homology.  
\end{abstract}

\maketitle


\section{Introduction}

In~\cite{khovanovfrob}, following his own seminal work in~\cite{khovanovsl2} 
and Lee, Bar-Natan and Turner's subsequent 
contributions~\cite{lee,bar-natancob,turner}, Khovanov classified all possible 
Frobenius systems of dimension two which give rise to 
link homologies via his construction in~\cite{khovanovsl2} and showed that 
there is a universal 
one, given by 
$$\bZ[X,a,b]/\left(X^2-aX-b\right).$$
Working over $\bC$, one can take $a$ and $b$ to be complex numbers and study 
the corresponding homology with coefficients in $\bC$. We refer to the latter 
as the $sl_2$-link homologies over $\bC$, because they are all 
related to the Lie algebra $sl_2$ (see \cite{khovanovfrob}).
Using the ideas in \cite{khovanovfrob,lee,turner}, it was shown 
in~\cite{mackaay-turner-vaz} that there are only two isomorphism 
classes of $sl_2$-link homologies over $\bC$. Given $a,b\in\bC$, 
the isomorphism class of 
the corresponding link homology is completely determined by the number of 
distinct roots of the polynomial $X^2-aX-b$. The original Khovanov $sl_2$-link homology $\mbox{KH}(L,\bC)$ corresponds to the choice $a=b=0$.

Bar-Natan~\cite{bar-natancob} obtained the universal $sl_2$-link homology 
in a different way, using a clever setup with cobordisms modulo relations. 
He shows how Khovanov's original construction of the $sl_2$-link 
homology~\cite{khovanovsl2} can be used to define a universal functor 
$\cU$ from the category of links, with link cobordisms modulo 
ambient isotopy as morphisms, to the homotopy category of complexes with 
values in the category of $1+1$-dimensional cobordisms modulo a finite set 
of universal relations. 
In the same paper 
he introduces the {\em tautological homology 
construction}, which produces an honest homology theory from $\mathcal U$. 
To obtain a finite dimensional homology one has to impose the 
extra relations
$$
\raisebox{-5pt}{
\includegraphics[height=0.2in]{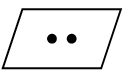}}=
a \raisebox{-5pt}{
\includegraphics[height=0.2in]{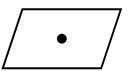}}+
b\raisebox{-5pt}{
\includegraphics[height=0.2in]{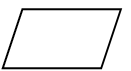}}
\quad\text{and}\quad
\raisebox{-17pt}{
\includegraphics[height=0.5in,width=0.2in]{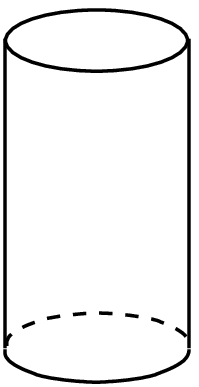}}=
\raisebox{-17pt}{
\includegraphics[height=0.5in,width=0.2in]{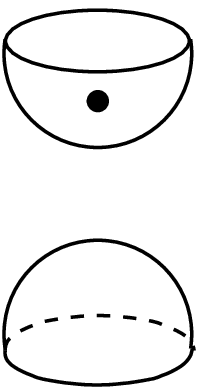}}+
 \raisebox{-17pt}{
\includegraphics[height=0.5in,width=0.2in]{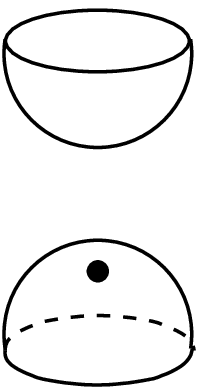}}
-a
\raisebox{-17pt}{
\includegraphics[height=0.5in,width=0.2in]{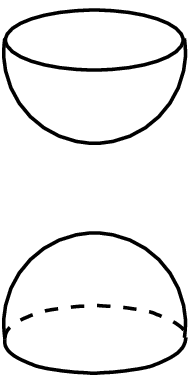}}
$$
on the cobordisms. 

In~\cite{khovanovsl3} Khovanov showed how to construct a link homology 
related to the Lie algebra $sl_3$. Instead of 
$1+1$-dimensional cobordisms, he uses webs and singular 
cobordisms modulo a finite 
set of relations, one of which is $X^3=0$. 
Gornik~\cite{gornik} studied the case when $X^3=1$, which is the 
analogue of Lee's work for $sl_3$. To be precise, Gornik studied a 
deformation of the Khovanov-Rozansky theory~\cite{KR} for $sl_n$, 
where $n$ is arbitrary. Khovanov and Rozansky followed a different 
approach to link homology using matrix factorizations which conjecturally 
yields the same for $sl_3$ as Khovanov's approach using 
singular cobordisms modulo relations~\cite{khovanovsl3}. 
However, in this paper we restrict to $n=3$ and only consider Gornik's 
results for this case.  

In the first part of this paper we construct the universal 
$sl_3$-link homology over $\bZ[a,b,c]$. For this universal construction 
we rely heavily on Bar-Natan's~\cite{bar-natancob} work on the universal 
$sl_2$-link homology and Khovanov's~\cite{khovanovsl3} work on his original 
$sl_3$-link homology. We first impose a finite set of relations 
on the category of webs and foams, analogous to 
Khovanov's~\cite{khovanovsl3} relations for his $sl_3$-link 
homology. We show that these relations imply certain identities between foams 
which are defined over $\bZ$ and are analogous to Bar-Natan's 
universal relations for $sl_2$ in \cite{bar-natancob}. The latter are 
sufficient to obtain a chain complex for each link which is invariant under 
the Reidemeister moves up to homotopy. However, they are insufficient 
for extending our construction to a functor, defined up to $\pm 1$, 
from the category of links to the homotopy category of complexes, 
for which we need the full set of 
relations including the ones depending on $a,b$ and $c$. 
\footnote{We thank M. Khovanov for spotting this problem in an earlier 
version of our paper.} This is the 
major difference with Bar-Natan's approach to the universal $sl_2$-homology. 
To obtain a finite-dimensional homology from our complex 
we use the tautological homology construction like 
Khovanov did in \cite{khovanovsl3} (the name {\em tautological homology} 
was coined by Bar-Natan in \cite{bar-natancob}). We denote this 
universal $sl_3$-homology by $U_{a,b,c}(L)$, 
which by the previous results is an invariant of the link $L$.  
 
In the second part of this paper we work over $\bC$ and take $a,b,c$ to 
be complex numbers, rather than formal parameters. We show that there 
are three isomorphism 
classes of $\uht * L$, depending on the number of 
distinct roots of the polynomial $f(X)=X^3-aX^2-bX-c$, and study them in 
detail. If $f(X)$ has only 
one root, with multiplicity three of course, then $\uht * L$ is isomorphic 
to Khovanov's original $sl_3$-link homology, which in our notation corresponds 
to $\uhtabc 0 0 0 * L$. If $f(X)$ has three distinct roots, then $\uht * L$ 
is isomorphic to Gornik's $sl_3$-link homology~\cite{gornik}, which corresponds 
to $\uhtabc 0 0 1 * L$. The case in which 
$f(X)$ has two distinct roots, one of which has multiplicity two, 
is new and had not been studied 
before to our knowledge, although in~\cite{DGR} and~\cite{GW} the authors 
make conjectures 
which are compatible with our results. We prove that there is a degree-preserving isomorphism
$$\uht * L\cong \bigoplus_{L'\subseteq L}\mbox{KH}^{*-j(L')}(L',\bC),$$
where $j(L')$ is a shift of degree $2\lk(L',L\backslash L')$. 
This isomorphism does not take into account the internal grading of the 
Khovanov homology.  

We have tried to make the paper reasonably self-contained, but we do assume 
familiarity with the papers by Bar-Natan~\cite{bar-natan,bar-natancob}, 
Gornik~\cite{gornik} and Khovanov~\cite{khovanovsl2,
khovanovsl3,khovanovfrob}.


\section{The universal $sl_3$-link homology}\label{sec:univhom}

\begin{figure}[h]
$$\xymatrix@C=16mm@R=2mm{
  &  \includegraphics[height=0.3in]{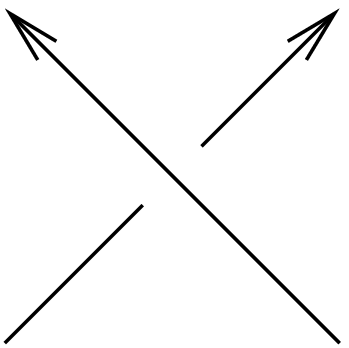} \ar[ld]_0\ar[rd]^1& \\
\includegraphics[height=0.3in]{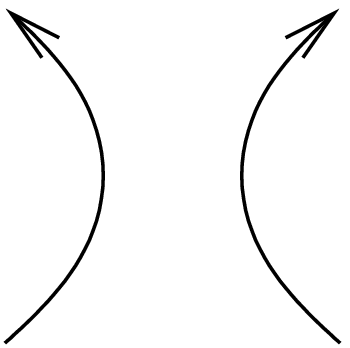} & & 
\includegraphics[height=0.3in]{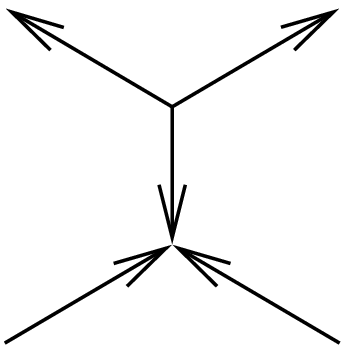} \\
  & \includegraphics[height=0.3in]{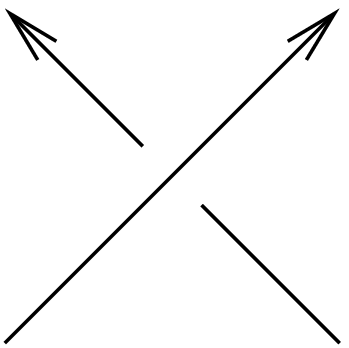} \ar[ul]^1 \ar[ur]_0&
}$$
\caption{0 and 1 resolutions of crossings}
\label{fig:resolutions}
\end{figure}

\begin{figure}[h]
\includegraphics[height=0.4in]{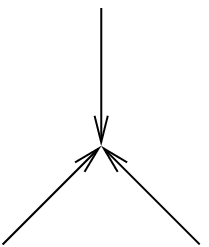}\qquad\quad
\includegraphics[height=0.4in]{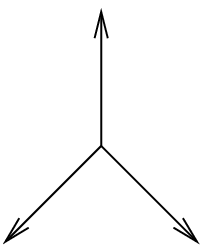}
\caption{``In'' and ``out'' orientations near a vertex}
\label{fig:in-out}
\end{figure}

\begin{figure}[h]
\includegraphics[height=0.4in]{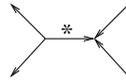}
\caption{Marked edges corresponding to a crossing in $D$}
\label{fig:marked-edge}
\end{figure}

\begin{figure}[h]
\includegraphics[height=0.5in]{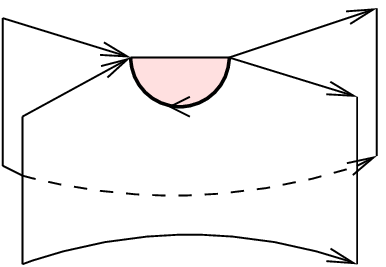}\qquad
\includegraphics[height=0.5in]{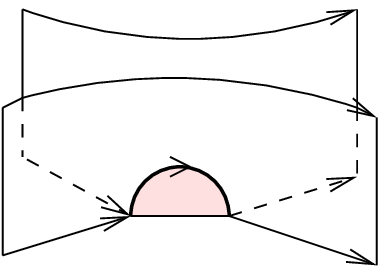}
\caption{Basic foams from $\orsmoothing$ to $\unorsmoothing$ (left) and from $\unorsmoothing$ to $\orsmoothing$ (right)}
\label{fig:ssaddle}
\end{figure}

Let $L$ be a link in $S^3$ and $D$ a diagram of $L$. In~\cite{khovanovsl3} Khovanov constructed a homological link invariant associated to $sl_3$. The construction starts by resolving each crossing of $D$ in two different ways, as in figure~\ref{fig:resolutions}. 
\noindent A diagram $\Gamma$ obtained by resolving all crossings of $D$ is an example of a \emph{web}. 
A web is a trivalent planar graph where near each vertex all the edges are oriented ``in'' or ``out'' 
(see figure~\ref{fig:in-out}). We also allow webs without vertices, which are 
oriented loops. 
Note that by definition our webs are closed; there are no vertices with fewer than 3 edges.  
Whenever it is necessary to keep track of crossings after their resolution we mark the corresponding edges as in figure~\ref{fig:marked-edge}.
A \emph{foam} is a cobordism with singular arcs between two webs. A singular arc in a foam $f$ is the set of points of $f$ that have a neighborhood homeomorphic to the letter Y times an interval (see the examples in figure~\ref{fig:ssaddle}).
Interpreted as morphisms, we read foams from bottom to top by convention, and the orientation of the singular arcs is by convention as in figure~\ref{fig:ssaddle}. Foams can have dots that can move freely on the facet to which they belong, but are not allowed to cross singular arcs. Let $\bZ[a,b,c]$ be the ring of 
polynomials in $a,b,c$ with integer coefficients. 

\begin{defn}
${\bf Foam}$ is the category whose objects are (closed) webs and whose morphisms are $\bZ[a,b,c]$-linear 
combinations of isotopy classes of foams.
\end{defn}

\noindent ${\bf Foam}$ is an additive category. For further details about this category 
see~\cite{khovanovsl3}.

From all different resolutions of all the crossings in $D$ we form a commutative 
hypercube of resolutions as in~\cite{khovanovsl3}. It has a web in each vertex and 
to an edge between two vertices, given by webs that differ only inside 
a disk $\mathcal D$ around one of the crossings of $D$, we associate the foam that is 
the identity everywhere except inside the cylinder ${\mathcal D}\times I$, where it 
looks like one of the basic foams in figure~\ref{fig:ssaddle}. An appropriate 
distribution of minus signs among the edges of the hypercube results in a chain 
complex of web diagrams analogous to the one in~\cite{bar-natancob} which we call 
$\brak{D}$, with ``column vectors'' of webs as ``chain objects'' and 
``matrices of foams'' as ``differentials''. We borrow some of the notation 
from~\cite{bar-natancob} and denote by $\mathbf{Kom}(\mathbf{Foam})$ the category of 
complexes in $\mathbf{Foam}$.

In subsections~\ref{ssec:loc_rel}-\ref{ssec:funct} we first impose a set of local relations on 
$\mathbf{Foam}$. We call this set of relations $\ell$ and denote by 
$\mathbf{Foam}_{/\ell}$ the category $\mathbf{Foam}$ divided by $\ell$. 
We prove a finite set of identities in $\mathbf{Foam}_{/\ell}$ 
which are universal in the sense that they are defined over $\bZ$. We then prove that the latter 
guarantee the invariance of 
$\brak{D}$ under the Reidemeister moves up to homotopy in $\mathbf{Kom}(\mathbf{Foam}_{/\ell})$ 
in a pictorial way, which is analogous to Bar-Natan's proof in~\cite{bar-natancob}. 
Note that the category $\mathbf{Kom}(\mathbf{Foam}_{/\ell})$ is analogous to Bar-Natan's category 
$Kob(\emptyset)=Kom(Mat(Cob^3_{/l}(\emptyset)))$. Next we show that up to signs 
$\brak{\;}$ is functorial under 
link cobordisms, i.e. defines a functor from $\mathbf{Link}$ to 
$\mathbf{Kom}_{/\pm h}(\mathbf {Foam}_{/\ell})$.
Here $\mathbf{Link}$ is the category of links in $S^3$ and ambient isotopy classes of link cobordisms 
properly embedded in $S^3\times[0,1]$ and $\mathbf{Kom}_{/\pm h}(\mathbf {Foam}_{/\ell})$ is the 
homotopy category of $\mathbf{Kom}(\mathbf {Foam}_{/\ell})$ modded out by $\pm 1$. 
For the functoriality we need all relations in $\ell$, including the ones which 
involve $a,b,$ and $c$. In subsection~\ref{ssec:univhom} we define a functor between 
$\mathbf{Foam}_{/\ell}$ and $\mathbf{\bZ[a,b,c]-{\text Mod}}$, the category of $\bZ[a,b,c]$-modules, 
which induces a homology functor 
$$U_{a,b,c}\colon \mathbf{Link}\ra\mathbf{Kom}_{/\pm h}(\mathbf{\bZ[a,b,c]-{\text Mod}}).$$

The principal ideas in this section, as well as most homotopies, are motivated by 
the ones in Khovanov's paper~\cite{khovanovsl3} and 
Bar-Natan's paper~\cite{bar-natancob}.

\subsection{Universal local relations}\label{ssec:loc_rel}
In order to construct the universal theory we divide $\mathbf{Foam}$ by the local relations  $\ell=(3D, CN, S, \Theta)$ below. 
$$\xymatrix@R=2mm
{
\raisebox{-5pt}{
\includegraphics[height=0.2in]{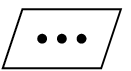}}=
a\raisebox{-5pt}{
\includegraphics[height=0.2in]{plan2dot}}+
b\raisebox{-5pt}{
\includegraphics[height=0.2in]{plan1dot}}+
c\raisebox{-5pt}{
\includegraphics[height=0.2in]{plan0dot}}
 & \text{(3D)}
\\
-\raisebox{-17pt}{
\includegraphics[height=0.5in,width=0.2in]{cylinder}}=
\raisebox{-17pt}{
\includegraphics[height=0.5in,width=0.2in]{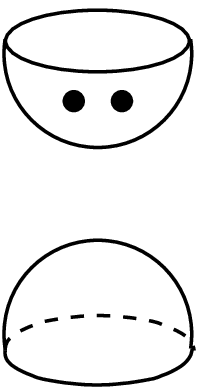}}+
\raisebox{-17pt}{
\includegraphics[height=0.5in,width=0.2in]{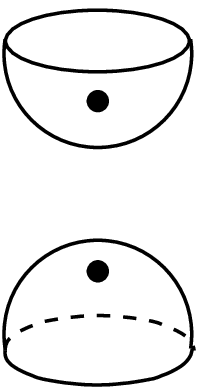}}+
\raisebox{-17pt}{
\includegraphics[height=0.5in,width=0.2in]{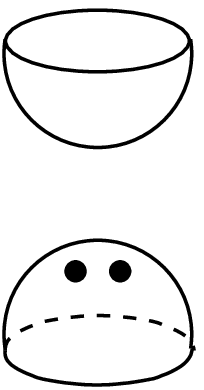}}
-a
\left( 
\raisebox{-17pt}{
\includegraphics[height=0.5in,width=0.2in]{cnecka1}}+
 \raisebox{-17pt}{
\includegraphics[height=0.5in,width=0.2in]{cnecka2}}
\right)
-b
\raisebox{-17pt}{
\includegraphics[height=0.5in,width=0.2in]{cneckb}}
 & \text{(CN)} \\
\raisebox{-8pt}{\includegraphics[width=0.3in]{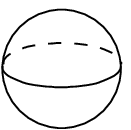}}=
\raisebox{-8pt}{\includegraphics[width=0.3in]{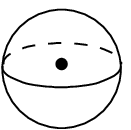}}=0,\quad
\raisebox{-8pt}{\includegraphics[width=0.3in]{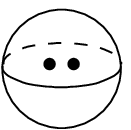}}=-1
 & \text{(S)}
}$$

\parpic[r]{
\includegraphics[height=0.4in]{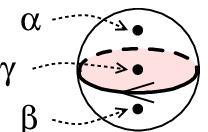}
} Recall from~\cite{khovanovsl3} that theta-foams are obtained by gluing three oriented disks along their boundaries (their orientations must coincide), as shown on the right. Note the orientation of the singular circle. Let $\alpha$, $\beta$, $\gamma$ denote the number of dots on each facet.
The ($\Theta$) \emph{relation} says that for $\alpha$, $\beta$ or $\gamma\leq 2$
\begin{equation*}
\theta(\alpha,\beta,\gamma)=\left\{
\begin{array}{cl}
1 & (\alpha,\beta,\gamma)=(1,2,0)\mbox{ or a cyclic permutation} \\ 
-1 & (\alpha,\beta,\gamma)=(2,1,0)\mbox{ or a cyclic permutation} \\ 
0 & \mbox{else}
\end{array}
\right.
\qquad(\Theta)
\end{equation*}
Reversing the orientation of the singular circle reverses the sign of $\theta(\alpha,\beta,\gamma)$. 
Note that when we have three or more dots on a facet of a foam we can use the (\emph{3D}) relation to reduce to the case where it has less than three dots.

A closed foam $f$ can be viewed as a morphism from the empty web to itself which by the relations  (\emph{3D}, \emph{CN}, \emph{S}, $\Theta$) is an element of $\bZ[a,b,c]$. It can be checked, as in ~\cite{khovanovsl3}, that this set of relations is consistent and determines uniquely the evaluation of every closed foam 
$f$, denoted $C(f)$. Define a $q$-grading on $\bZ[a,b,c]$ as $q(1)=0$, $q(a)=2$, $q(b)=4$ and $q(c)=6$. As in~\cite{khovanovsl3} we define the \emph{q-grading} of a foam $f$ with $d$ dots by $$q(f)=-2\chi(f)+ \chi(\partial f)+2d,$$ where $\chi$ denotes the Euler characteristic.

\begin{defn}
${\bf Foam}_{/\ell}$ is the quotient of the category ${\bf Foam}$ by the local relations $\ell$. 
For webs $\Gamma$, $\Gamma'$ and for families $f_i\in\Hom_{\mathbf{Foam}_{/\ell}}(\Gamma,\Gamma')$ and 
$c_i\in\bZ[a,b,c]$ we impose $\sum_ic_if_i=0$ if and only if 
$\sum_ic_iC(g'f_ig)=0$ holds, for all 
$g\in\Hom_{\mathbf{Foam}_{/\ell}}(\emptyset,\Gamma)$ and 
$g'\in\Hom_{\mathbf{Foam}_{/\ell}}(\Gamma',\emptyset)$.
\end{defn}

\begin{lem}\label{lem:identities}
We have the following relations in ${\bf Foam}_{/\ell}$:
$$\xymatrix@R=2mm{
\raisebox{-16pt}{\includegraphics[width=0.6in]{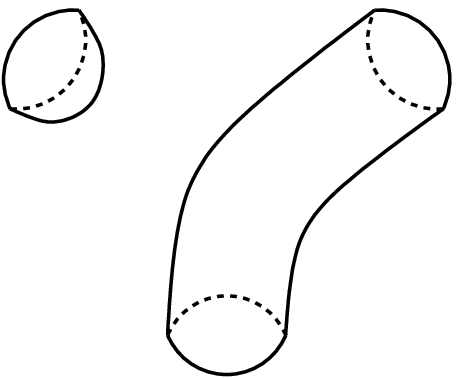}}+
\raisebox{-16pt}{\includegraphics[width=0.6in]{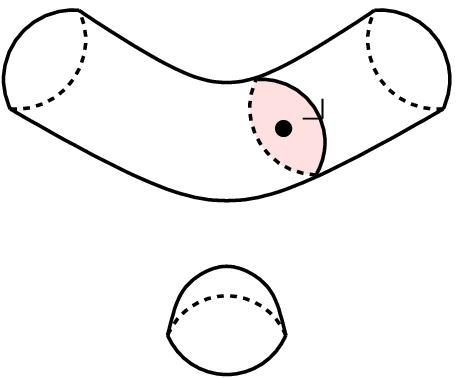}}=
\raisebox{-16pt}{\includegraphics[width=0.6in]{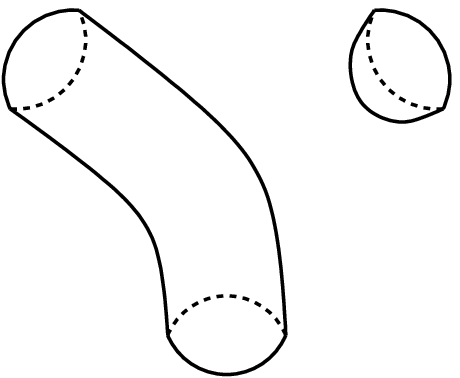}}+
\raisebox{-16pt}{\includegraphics[width=0.6in]{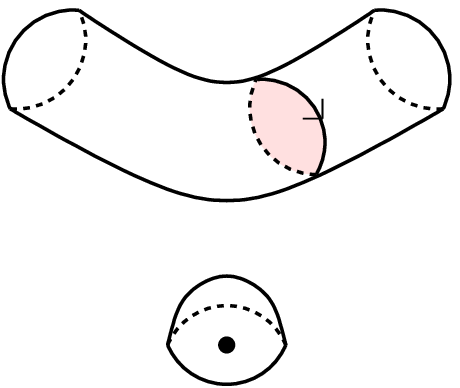}}
 & \text{(4C)}  \\
\label{R-rel}
\raisebox{-16pt}{\includegraphics[height=0.5in,width=0.2in]{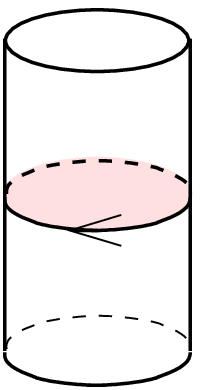}}=
\raisebox{-16pt}{\includegraphics[height=0.5in,width=0.2in]{cnecka1}}-
\raisebox{-16pt}{\includegraphics[height=0.5in,width=0.2in]{cnecka2}}
 & \text{(RD)}  \\
\raisebox{-20pt}{\includegraphics[height=0.6in]{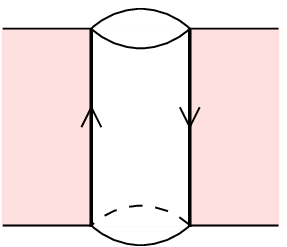}}=
\raisebox{-26pt}{\includegraphics[height=0.75in]{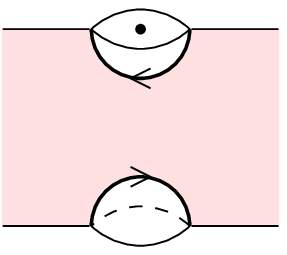}}-
\raisebox{-26pt}{\includegraphics[height=0.75in]{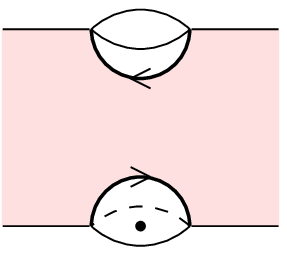}}
 & \text{(DR)}  \\
\raisebox{-28pt}{\includegraphics[height=0.8in]{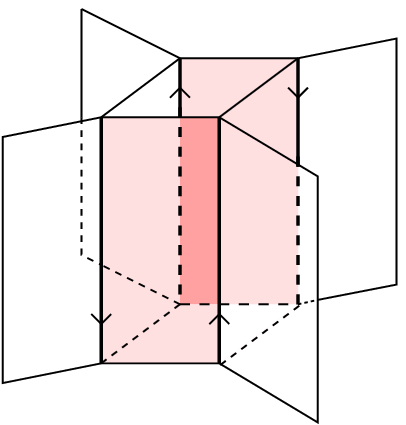}}= - 
\raisebox{-28pt}{\includegraphics[height=0.8in]{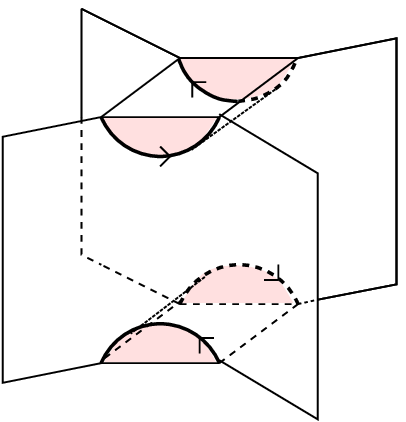}}- 
\raisebox{-28pt}{\includegraphics[height=0.8in]{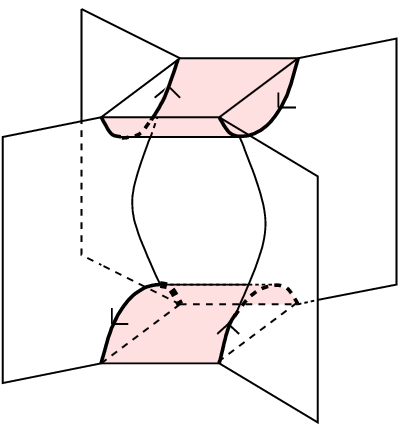}}
 & \text{(SqR)}
}$$
\end{lem}

\begin{proof}
Relations (\emph{4C}) and (\emph{RD}) are immediate and follow from $(\emph{CN})$ 
and ($\Theta$). Relations (\emph{DR}) and (\emph{SqR}) are proved as 
in~\cite{khovanovsl3} (see also lemma~\ref{lem:KhK})
\end{proof}

The following equality, and similar versions, which corresponds to an isotopy, 
we will often use in the sequel 
\begin{equation}\label{eq:usefuleq}
\raisebox{-0.3in}{\includegraphics[height=0.7in]{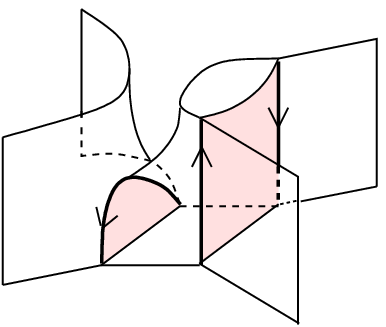}}
\circ
\raisebox{-0.3in}{\includegraphics[height=0.7in]{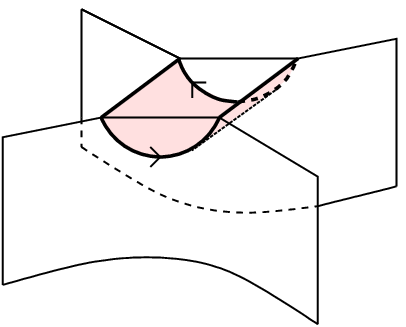}}
=
\raisebox{-0.3in}{\includegraphics[height=0.7in]{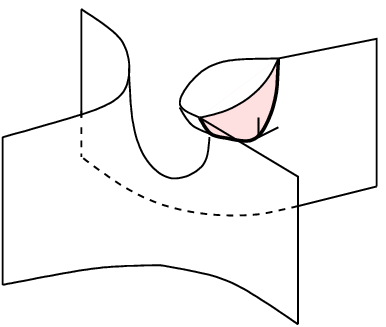}}
\end{equation}
where $\circ$ denotes composition of foams.

\begin{figure}[h]
$$\xymatrix@R=1mm{
\raisebox{-0.27in}{\includegraphics[height=0.5in]{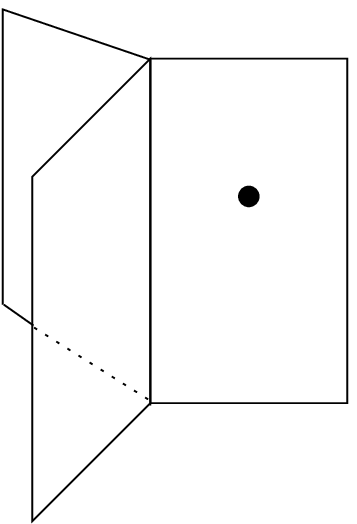}}
\,+\,
\raisebox{-0.27in}{\includegraphics[height=0.5in]{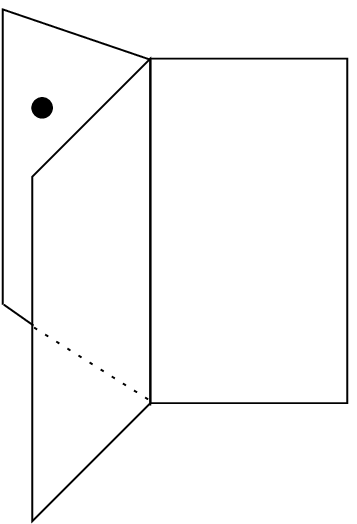}}
\,+\,
\raisebox{-0.27in}{\includegraphics[height=0.5in]{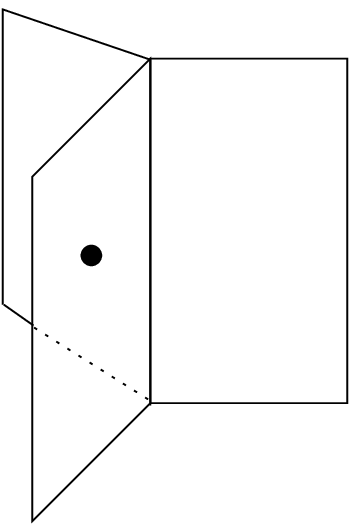}}
\,=\,a\,
\raisebox{-0.27in}{\includegraphics[height=0.5in]{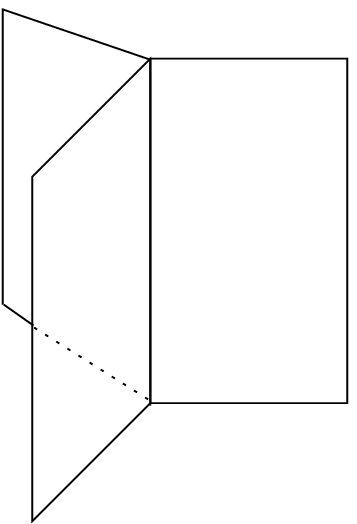}} 
\\
\raisebox{-0.27in}{\includegraphics[height=0.5in]{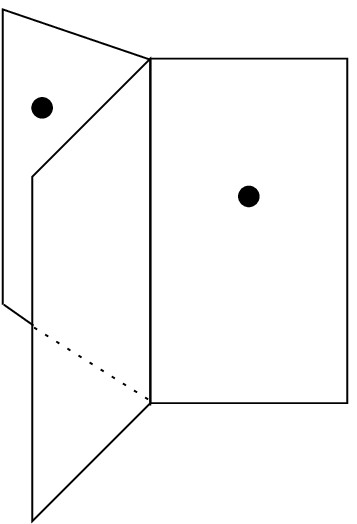}}
\,+\,
\raisebox{-0.27in}{\includegraphics[height=0.5in]{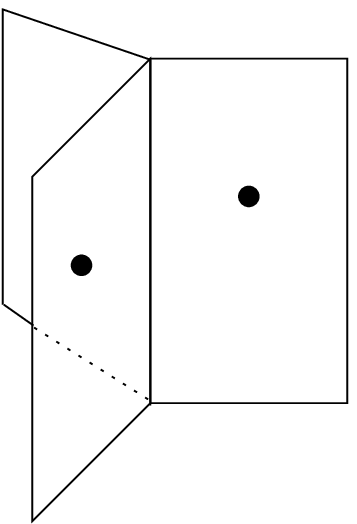}}
\,+\,
\raisebox{-0.27in}{\includegraphics[height=0.5in]{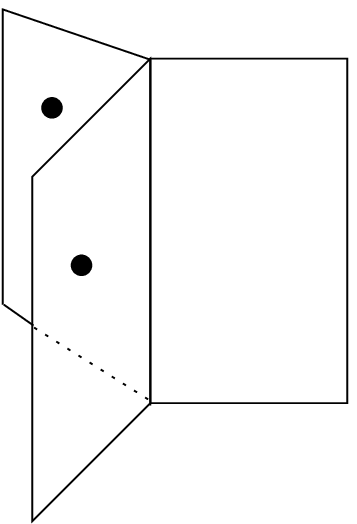}}
\,=\,-b\,
\raisebox{-0.27in}{\includegraphics[height=0.5in]{pdots000}} 
\\
\raisebox{-0.27in}{\includegraphics[height=0.5in]{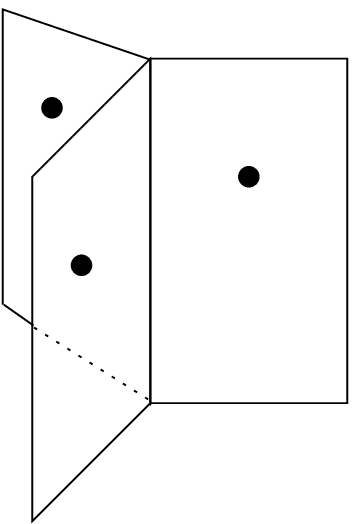}}
\,=\,c\,
\raisebox{-0.27in}{\includegraphics[height=0.5in]{pdots000}}
}$$
\caption{Exchanging dots between faces. The relations are the same regardless  of which edges are marked and the orientation on the singular arcs.}
\label{fig:pdots}
\end{figure}

In figure~\ref{fig:pdots} we also have a set of useful identities which establish 
the way we can exchange dots between faces. These identities can be used for the simplification of foams and are an immediate consequence of the relations in $\ell$.

\subsection{Invariance under the Reidemeister moves}\label{ssec:thm-inv}

In this subsection we prove invariance of $\brak{\;}$ under the Reidemeister moves. The proof only 
uses the relations which are defined over $\bZ$. The main result of this section is the following

\begin{thm}
$\brak{D}$ is invariant under the Reidemeister moves up to homotopy, in other words it is an invariant in $\mathbf{Kom}_{/h}(\mathbf{Foam}_{/\ell})$.
\label{thm-inv}
\end{thm}

\begin{proof}
To prove invariance under the Reidemeister moves we work diagrammatically and only 
use the identities in lemma~\ref{lem:identities} along with the (\emph{S}) relation.

\subsubsection*{Reidemeister I} 
Consider diagrams $D$ and $D'$ that differ only in a circular region as in the 
figure below.
$$D=\raisebox{-13pt}{\includegraphics[height=0.4in]{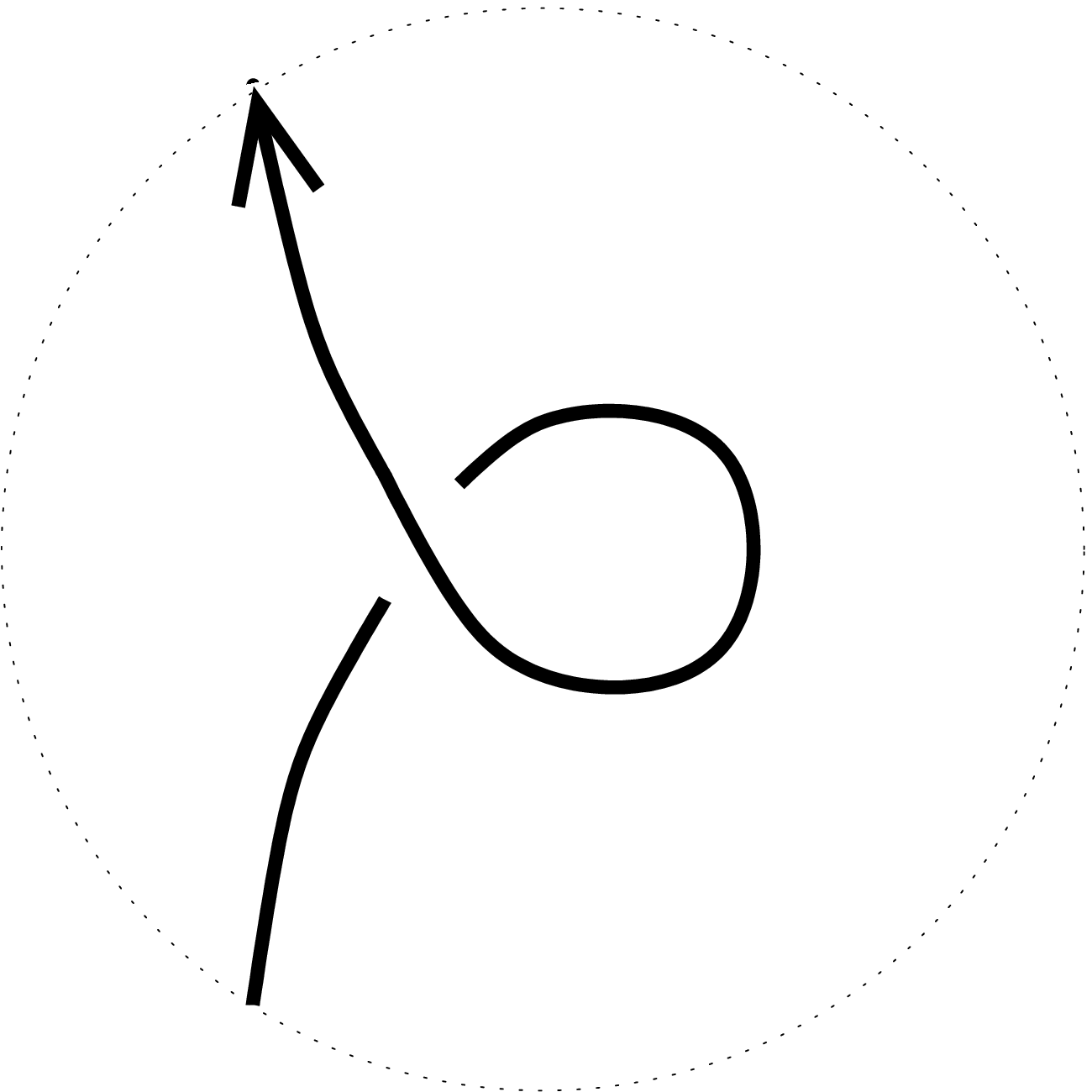}}\qquad
D'=\raisebox{-13pt}{\includegraphics[height=0.4in]{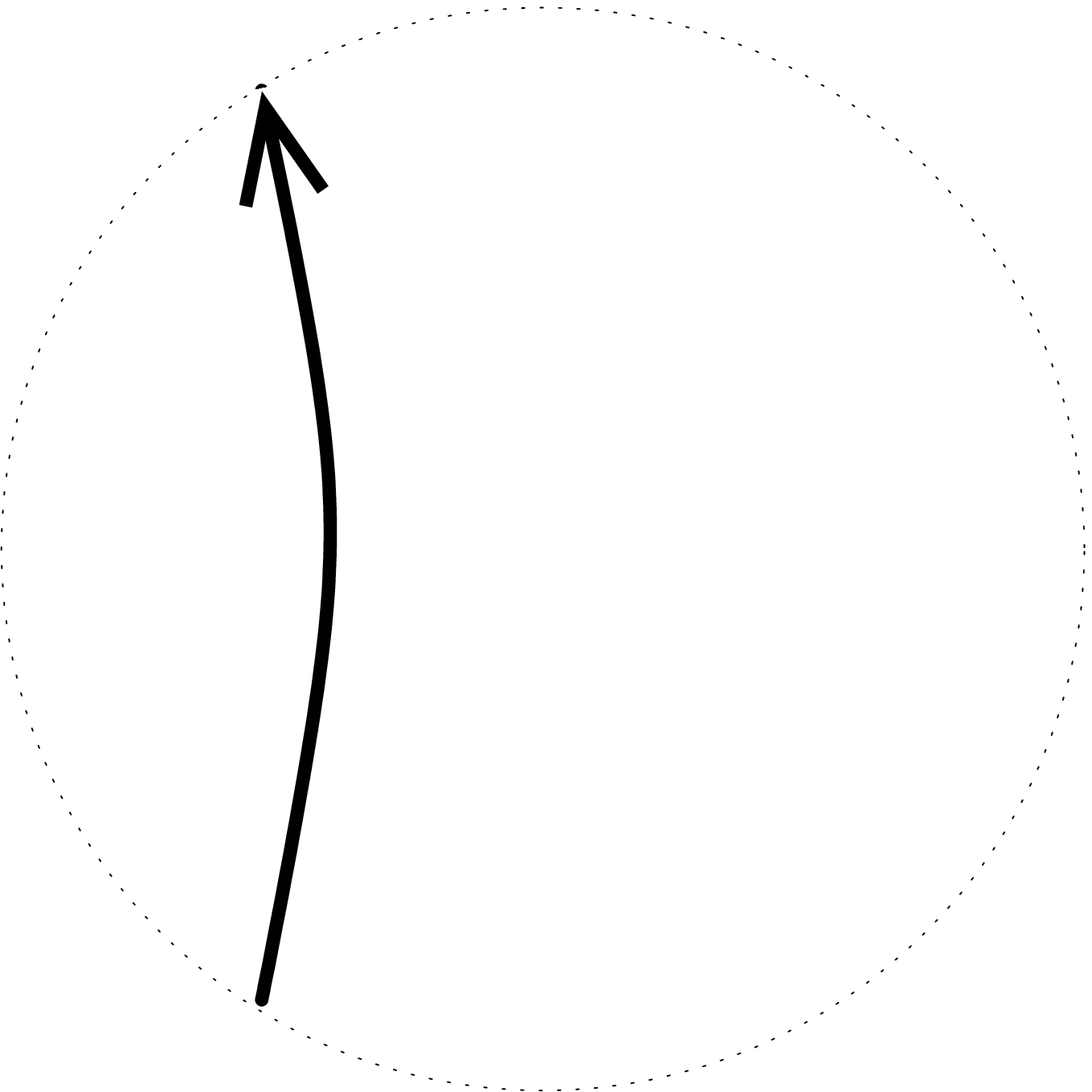}}\,$$
\begin{figure}[h]
$\xymatrix@R=32mm{
  \brak{D}:
\\
  \brak{D'}:
}
\xymatrix@C=35mm@R=28mm{
 \includegraphics[height=0.3in]{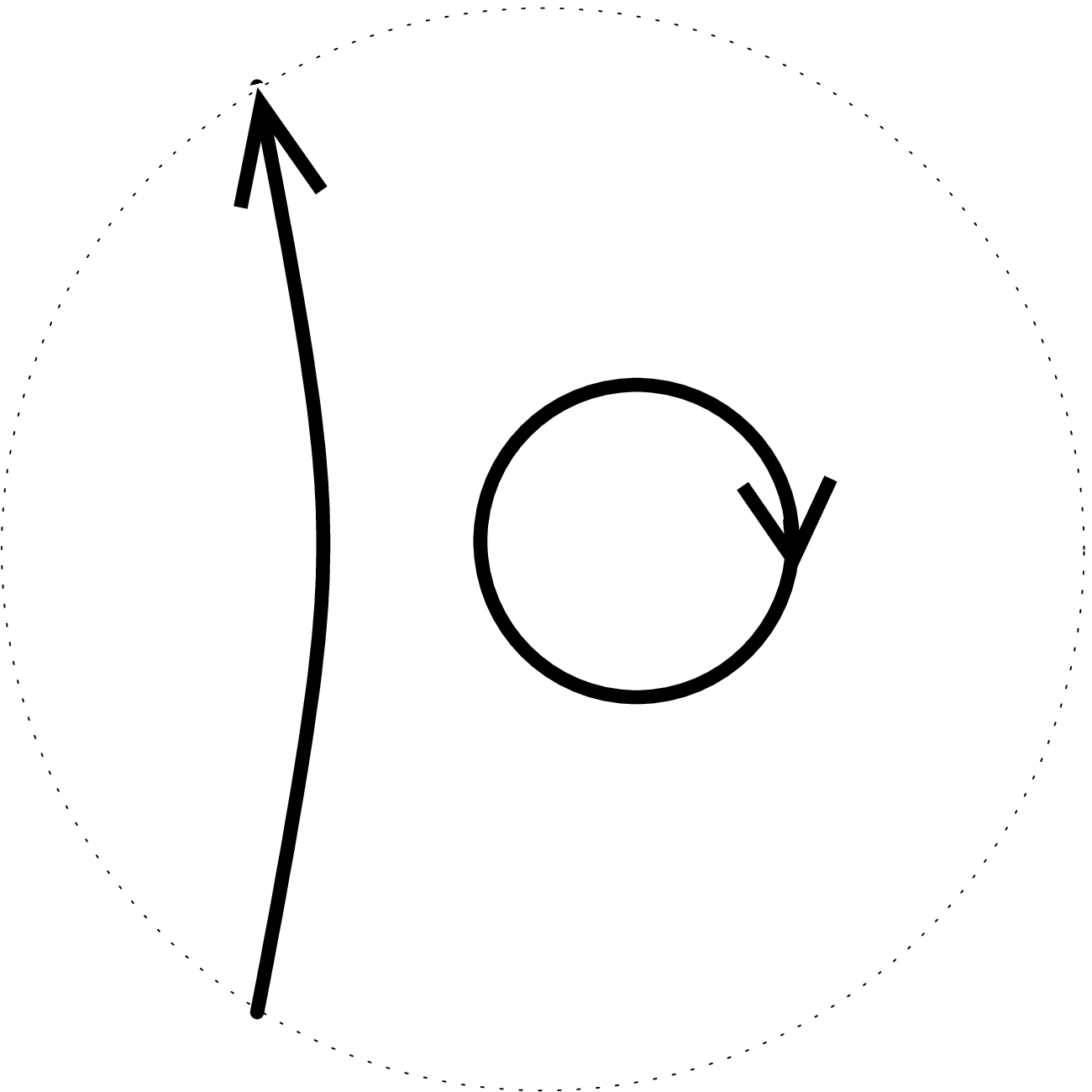}
  \ar@<4pt>[d]^{
        g^0 \;=\; \raisebox{-22pt}{\includegraphics[height=0.6in]{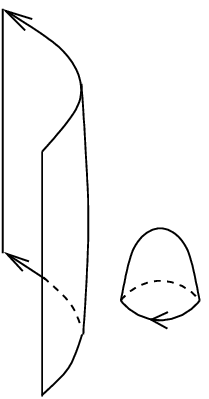}}} 
  \ar@<4pt>[r]^{
        d \;=\; \raisebox{-24pt}{\includegraphics[height=0.6in]{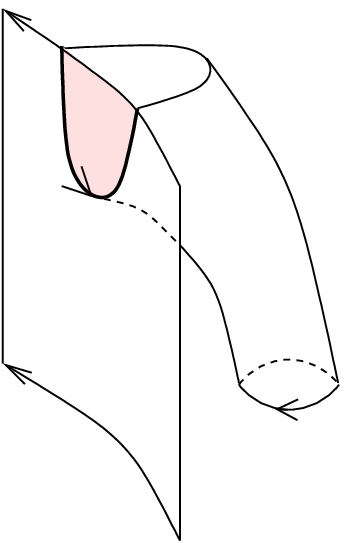}} } &
  \includegraphics[height=0.3in]{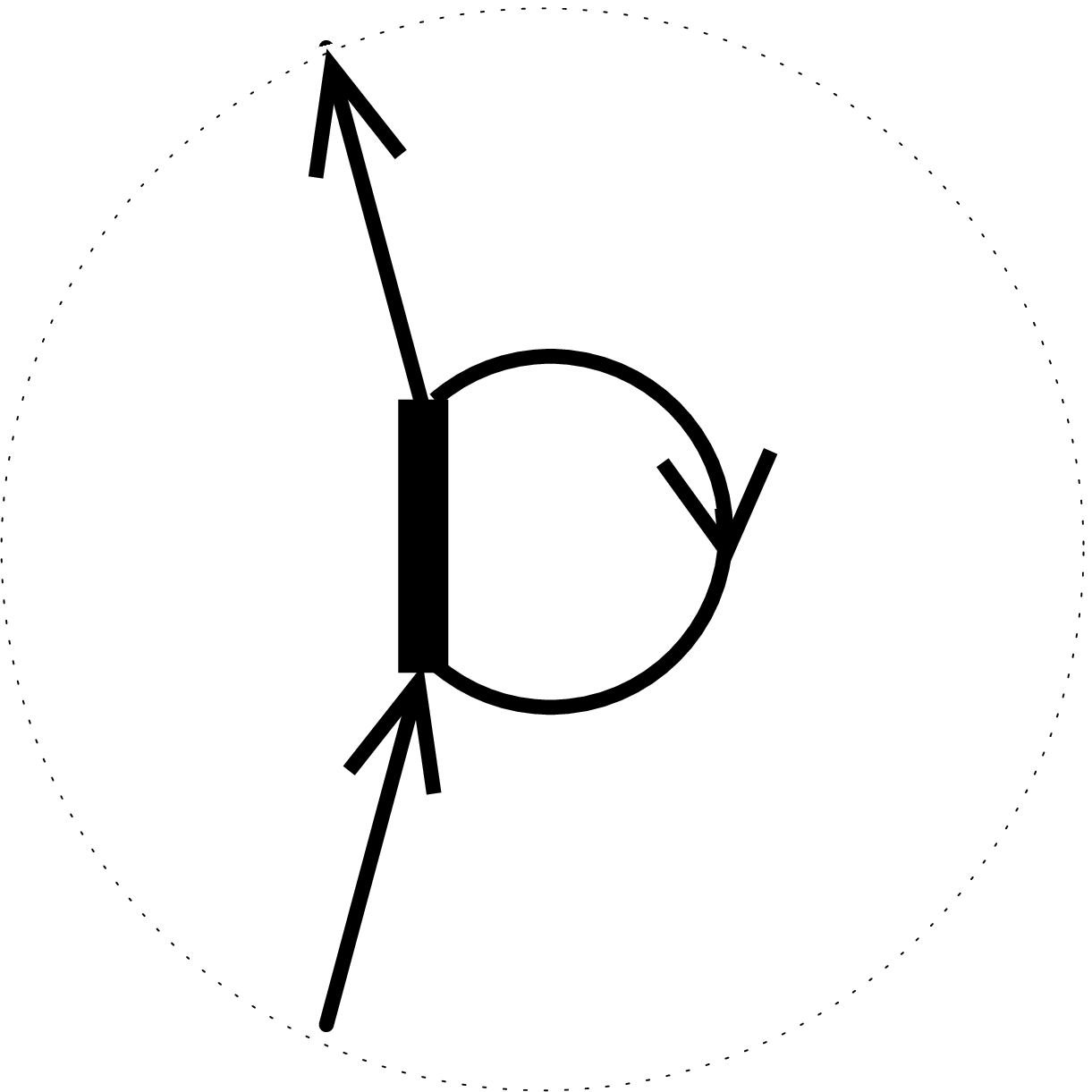} \ar@<2pt>[l]^{
     h \;=\; - \;\raisebox{-24pt}{\includegraphics[height=0.6in]{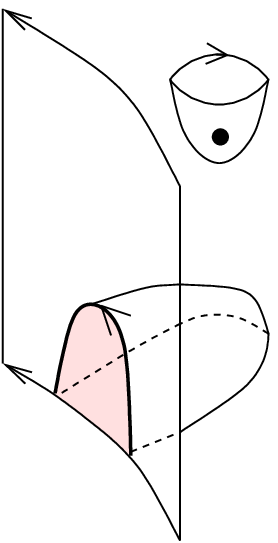}} \; + \;
            \raisebox{-24pt}{\includegraphics[height=0.6in]{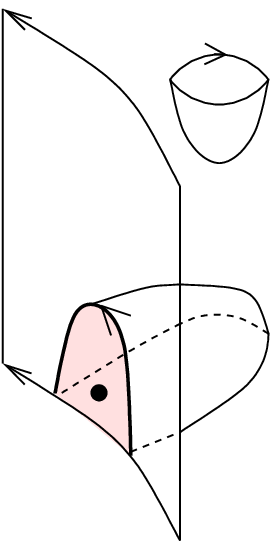}} }  \\
 \includegraphics[height=0.3in]{ReidI-1}\ar[r]^0 
  \ar@<4pt>[u]^{
        f^0 \;=\; \raisebox{-22pt}{\includegraphics[height=0.6in]{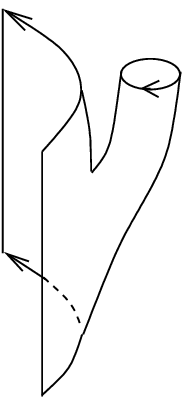}} } &
  0 \ar@{<->}[u]_0 
}$
\caption{Invariance under $Reidemeister\,I$} \label{fig:RI_Invariance}
\end{figure} 
\noindent We give the homotopy between complexes $\brak{D}$ and $\brak{D'}$ in figure~\ref{fig:RI_Invariance}. It is immediate that $g^0f^0=Id(\edge)$.
To see that $df^0=0$ use (\emph{DR}) near the top of $df^0$ and then (\emph{RD}).
The equality $dh=id(\digon)$ follows from (\emph{DR}) (note the orientations on the singular circles) and $f^0g^0+hd=id(\edgewcircle)$ follows from (\emph{4C}). Therefore $\brak{D'}$  is homotopy-equivalent to $\brak{D}$.

\subsubsection*{Reidemeister IIa}
Consider diagrams $D$ and $D'$ that differ in a circular region, as
in the figure below.
$$D=\raisebox{-13pt}{\includegraphics[height=0.4in]{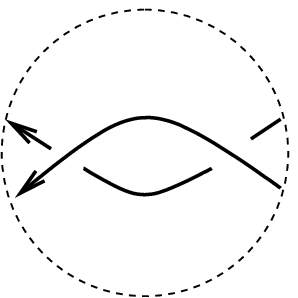}}\qquad
D'=\raisebox{-13pt}{\includegraphics[height=0.4in]{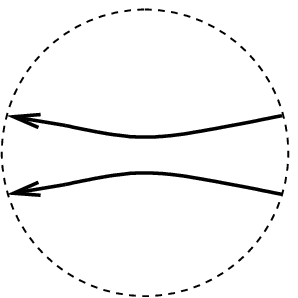}}$$
\begin{figure}[b]
\includegraphics[width=4.9in]{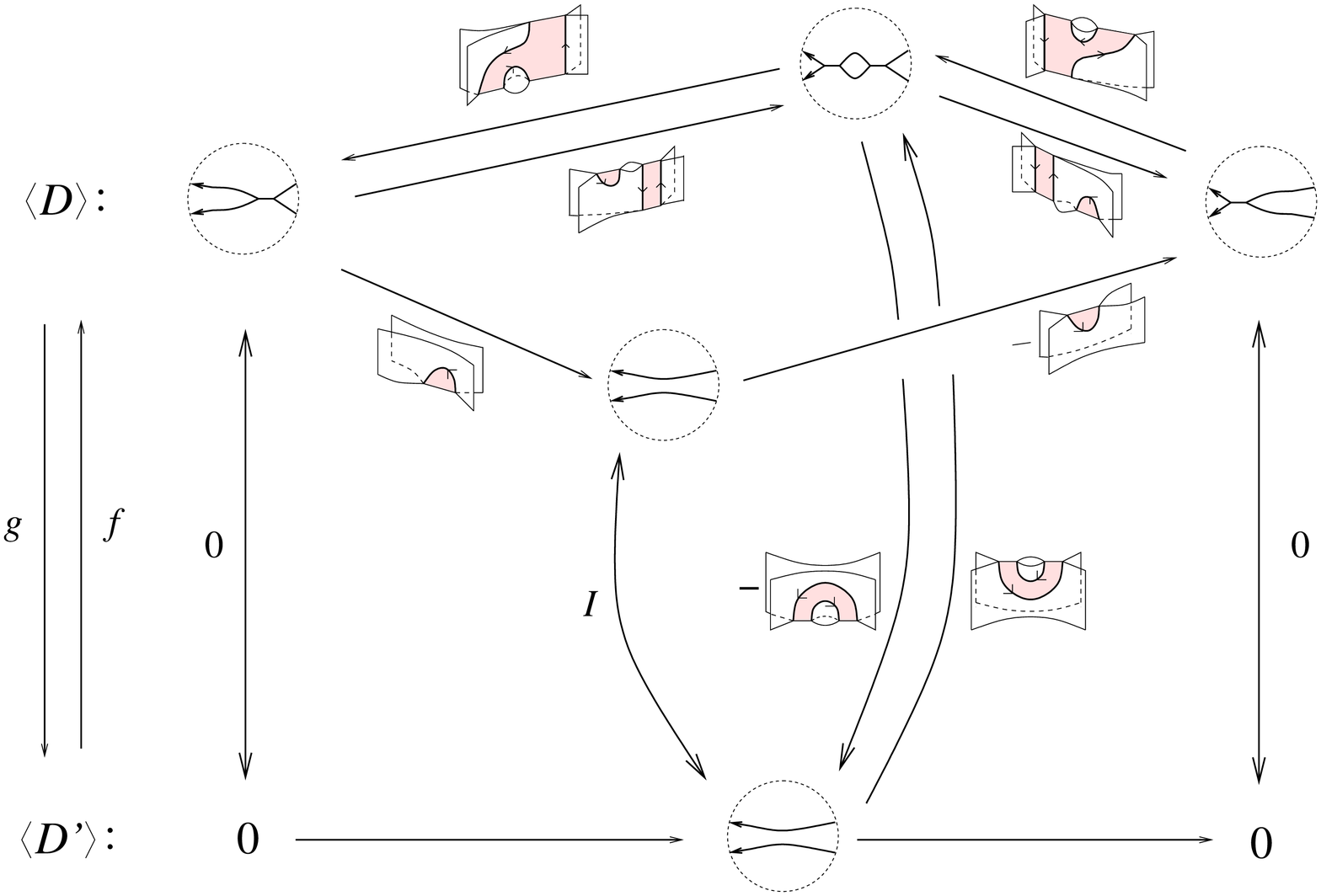}
\caption{Invariance under $Reidemeister\,IIa$} \label{fig:RIIaInvariance}
\end{figure}

We leave to the reader the task of checking that the diagram in 
figure~\ref{fig:RIIaInvariance} defines a homotopy between the 
complexes $\brak{D}$ and $\brak{D'}$: 
\begin{itemize}
\item $g$ and $f$ are morphisms of complexes (use only isotopies); 
\item $g^1f^1=Id_{\brak{D'}^1}$  (use (\emph{RD})); 
\item $f^0g^0+hd=Id_{\brak{D}^0}$ and $f^2g^2+dh=Id_{\brak{D}^2}$  (use isotopies);
\item $f^1g^1+dh+hd=Id_{\brak{D}^1}$  (use (\emph{DR})).
\end{itemize}

\subsubsection*{Reidemeister IIb}
Consider diagrams $D$ and $D'$ that differ only in a circular region, as 
in the figure below.
$$D=\raisebox{-13pt}{\includegraphics[height=0.4in]{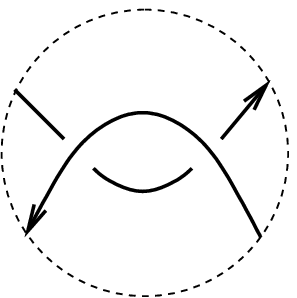}}\qquad
D'=\raisebox{-13pt}{\includegraphics[height=0.4in]{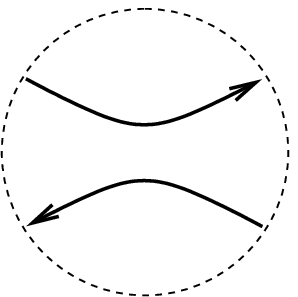}}$$
Again, checking that the diagram in figure~\ref{fig:RIIbInvariance} defines a 
homotopy between the complexes $\brak{D}$ and $\brak{D'}$ is left to the reader:
\begin{itemize}
\item $g$ and $f$ are morphisms of complexes (use only isotopies); 
\item $g^1f^1=Id_{\brak{D'}^1}$  (use (\emph{RD}) and (\emph{S})); 
\item $f^0g^0+hd=Id_{\brak{D}^0}$ and $f^2g^2+dh=Id_{\brak{D}^2}$ (use (\emph{RD}) and (\emph{DR}));
\item $f^1g^1+dh+hd=Id_{\brak{D}^1}$ (use (\emph{DR}), (\emph{RD}), (\emph{4C}) and (\emph{SqR})).
\end{itemize}
\begin{figure}[h]
\includegraphics[width=4.9in]{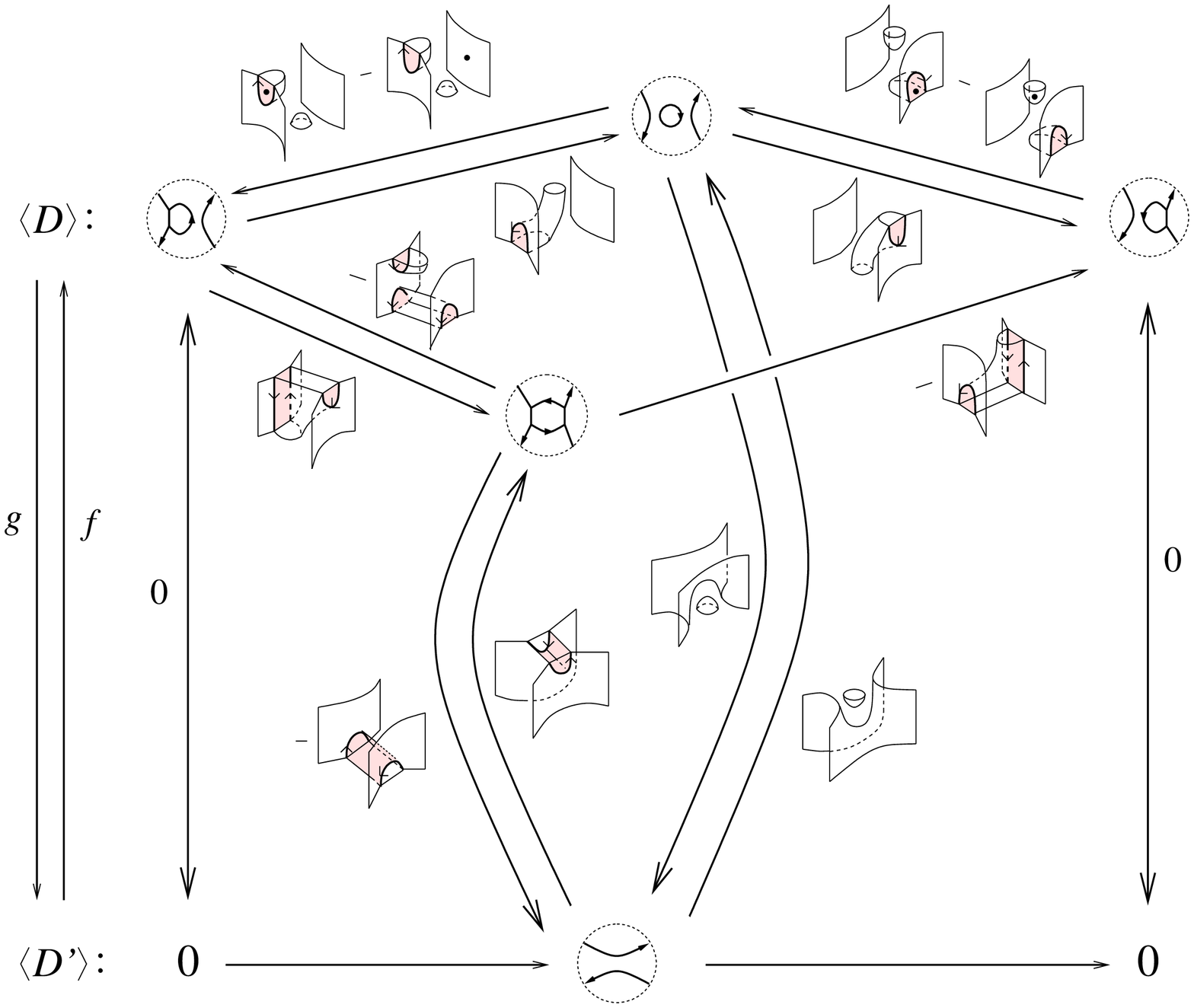}
\caption{Invariance under $Reidemeister\,IIb$} \label{fig:RIIbInvariance}
\end{figure}

\subsubsection*{Reidemeister III}
Consider diagrams $D$ and $D'$ that differ only in a circular region, as
in the figure below.
$$D=\raisebox{-13pt}{\includegraphics[height=0.4in]{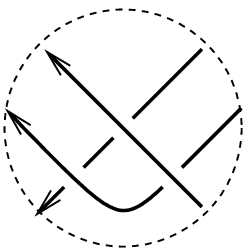}}\qquad
D'=\raisebox{-13pt}{\includegraphics[height=0.4in]{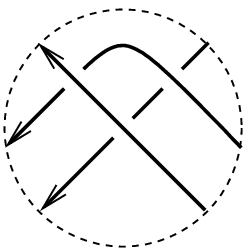}}$$
\begin{figure}[h]
\includegraphics[width=4.1in]{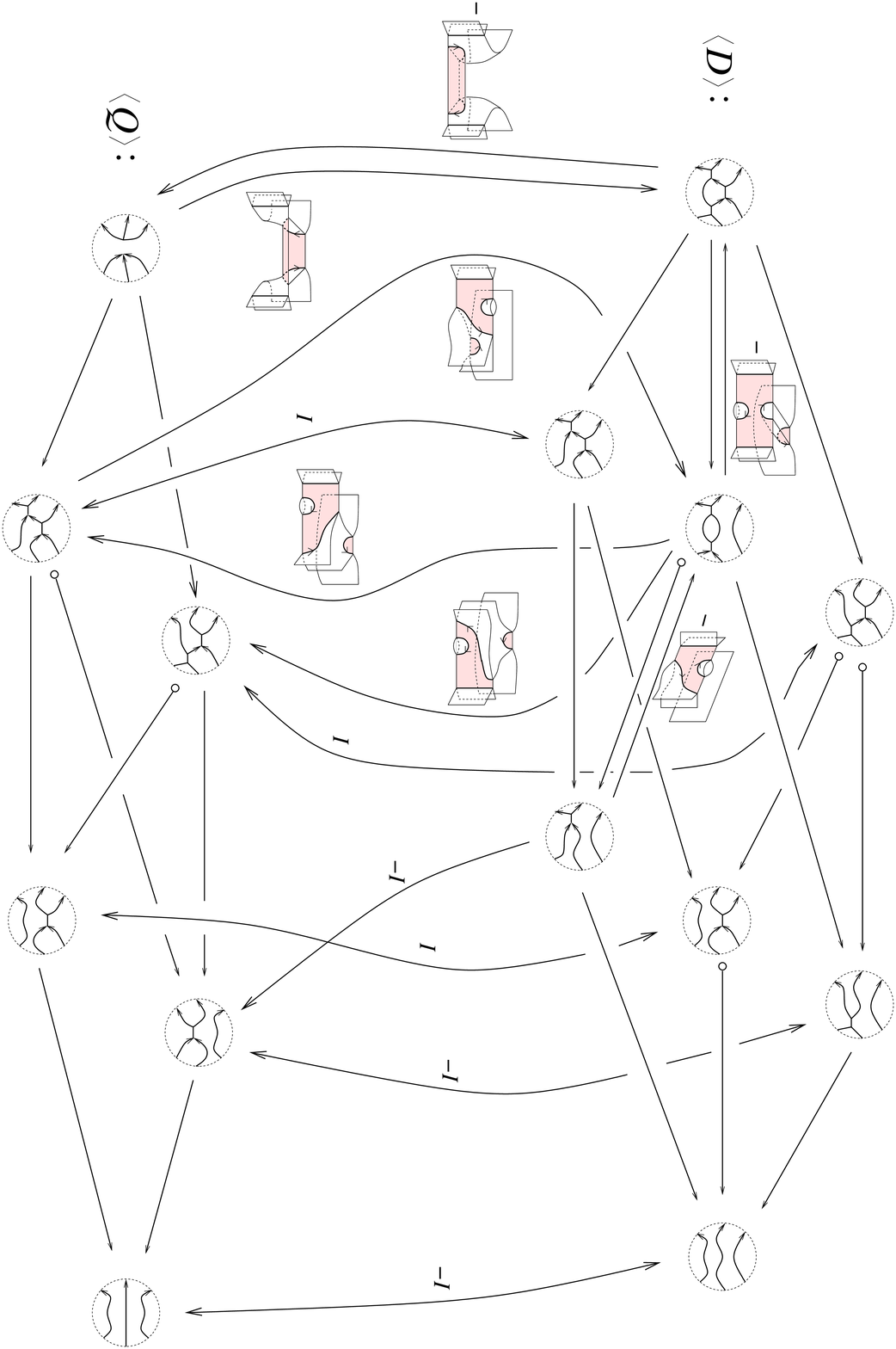}
\caption{First step of invariance under $Reidemeister\,III$. A circle attached to the tail of an arrow indicates that the corresponding morphism has a minus sign.}
\label{fig:RIIIInvariance}
\end{figure}

We prove that $\brak{D'}$ is homotopy equivalent to $\brak{D}$ by showing that both complexes are homotopy equivalent to a third complex denoted $\brak{Q}$ (the bottom complex in figure~\ref{fig:RIIIInvariance}).
Figure~\ref{fig:RIIIInvariance} shows that $\brak{D}$ is homotopy equivalent to $\brak{Q}$. By applying a symmetry relative to a horizontal axis crossing each diagram in $\brak{D}$ we obtain a homotopy equivalence between $\brak{D'}$ and $\brak{Q}$. It follows that $\brak{D}$ is homotopy equivalent to $\brak{D'}$.
\end{proof}

Theorem~\ref{thm-inv} allows us to use any diagram $D$ of $L$ to obtain the invariant in $\mathbf{Kom}_{/h}(\mathbf{Foam}_{/\ell})$ and justifies the notation $\brak{L}$ for $\brak{D}$.

\subsection{Functoriality}\label{ssec:funct}
It is clear that the construction and the results of the previous sections can be 
extended to the category of tangles, following Bar-Natan's approach 
in~\cite{bar-natancob}. One can then prove functoriality of $\brak{\;}$ as Bar-Natan does. 
The proofs of lemmas~8.6-8.8 in~\cite{bar-natancob} are identical. The proof of 
lemma~8.9 follows the same reasoning but uses the homotopies of our subsection~\ref{ssec:thm-inv}. Without giving any details of this generalization, we state the main result. Let $\mathbf{Kom}_{/\pm h}(\mathbf{Foam}_{/\ell})$ denote the category $\mathbf{Kom}_{/h}(\mathbf{Foam}_{/\ell})$ modded out by $\pm 1$. Then

\begin{prop}
\label{prop:func}
$\brak{\;}$ defines a functor $\mathbf{Link}\ra \mathbf{Kom}_{/\pm h}(\mathbf{Foam}_{/\ell})$.
\end{prop}

\subsection{Universal homology}\label{ssec:univhom}

Following Khovanov~\cite{khovanovsl3}, we now define a functor $C$ between $\mathbf{Foam}_{/\ell}$ and 
$\mathbf{\bZ[a,b,c]-{\text Mod}}$, which extends in a straightforward manner to the category $\mathbf{Kom}(\mathbf{Foam}_{/\ell})$.

\begin{defn}
For a closed web $\Gamma$, define $C(\Gamma)=\Hom_{\mathbf{Foam}_{/\ell}}(\emptyset,\Gamma)$. 
From the $q$-grading formula for foams, it follows that $C(\Gamma)$ is graded. 
For a foam $f$ between webs $\Gamma$ and $\Gamma'$ we define the $\bZ[a,b,c]$-linear map 
$$C(f):\Hom_{\mathbf{Foam}_{/\ell}}(\emptyset,\Gamma) \ra 
\Hom_{\mathbf{Foam}_{/\ell}}(\emptyset,\Gamma')$$ given by composition, whose degree equals $q(f)$. 
\end{defn}
Note that, if we have a disjoint union of webs $\Gamma$ and $\Gamma'$, then $C(\Gamma \sqcup \Gamma')\cong C(\Gamma)\otimes C(\Gamma')$.

The following relations are a categorified version of Kuperberg's skein relations~\cite{Kup} and were used and proved by Khovanov in~\cite{khovanovsl3} to relate his 
$sl_3$-link homology to the quantum $sl_3$-link invariant. 

\begin{lem}\label{lem:KhK}
{\em (Khovanov-Kuperberg relations~\cite{khovanovsl3, Kup})} We have the following decompositions under 
the functor $C$:
$$\xymatrix@R=3.8mm{
C(\unknot\Gamma) \cong
C(\unknot)\otimes C(\Gamma) &
\mbox{\emph{(Circle Removal)}}
\\
C(\raisebox{-2.6pt}{\includegraphics[height=0.15in]{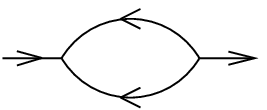}}) \cong 
C(\raisebox{2pt}{\includegraphics[height=0.025in,width=0.3in]{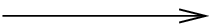}})\{-1\}\oplus C(\raisebox{2pt}{\includegraphics[height=0.025in,width=0.3in]{hthickedge}})\{1\} & \mbox{\emph{(Digon Removal)}}
\\
C\left(\raisebox{-8.0pt}{\includegraphics[height=0.3in]{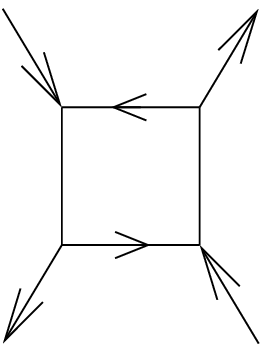}}\right)\cong 
C\left(\raisebox{-8.0pt}{\includegraphics[height=0.3in]{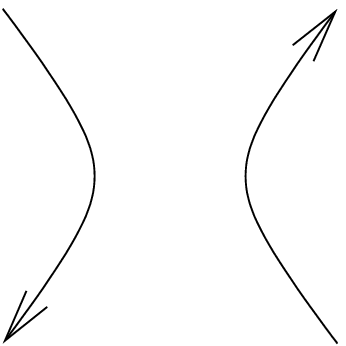}}\right)\oplus
C\left(\raisebox{-8.0pt}{\includegraphics[height=0.3in]{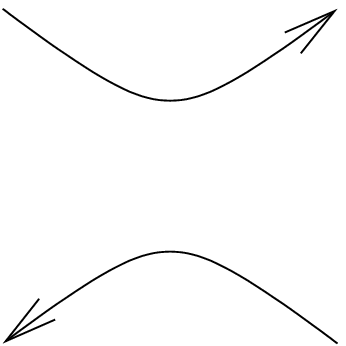}}\right) &
\mbox{\emph{(Square Removal)}}
}$$
where $\{j\}$ denotes a positive shift in the $q$-grading by $j$.
\end{lem}

\begin{proof}
\emph{Circle removal} is immediate from the definition of $C(\Gamma)$. \emph{Digon removal} and \emph{square removal} are proved as in~\cite{khovanovsl3}. Notice that \emph{Digon removal} and \emph{square removal} are related to the local relations (\emph{DR}) and (\emph{SqR}) of page~\pageref{R-rel}.
\end{proof}

Let $\mathcal H$ be the homology functor. We denote by $\uhtgen * D$ the composite functor 
${\mathcal H}^*C\brak{D}$. Proposition~\ref{prop:func} implies

\begin{prop}\label{prop:inv_vec}
$U_{a,b,c}:\mathbf{Link}\ra\mathbf{Kom_{/\pm h}(\bZ[a,b,c]-Mod)}$ defines a functor.
\end{prop}

\noindent We use the notation $C(L)$ for $C\brak{D}$ and $\uhtgen * L$ for $\uhtgen * D$.


\section{Isomorphism classes}\label{sec:isos}

In this section we work over $\bC$ and take $a,b,c$ to be complex numbers. 
Using the same construction as in the first part of this paper we can define $\uht * L$, which is 
the universal $sl_3$-homology with coefficients in $\bC$.
We show that there are three isomorphism classes of 
$\uht * L$. Throughout this section we write $f(X)=X^3-aX^2-bX-c$. For 
a given choice of $a,b,c\in\bC$, the isomorphism class of $\uht * L$ 
is determined by the number of distinct roots of $f(X)$. 

\begin{rem}
We could work over $\bQ$ just as well and obtain the same results, 
except that in the proofs we would first have to pass to quadratic 
or cubic field extensions of $\bQ$ to guarantee the existence of the 
roots of $f(X)$ in the field of coefficients of the homology. 
The arguments we present for $\uht * L$ remain valid over those quadratic 
or cubic extensions. The 
universal coefficient theorem then shows that our results hold true 
for the homology defined over $\bQ$.    
\end{rem}

If $f(X)=(X-\alpha)^3$, 
then the isomorphism $\planedot\mapsto\planedot - \alpha\plane$ induces 
an isomorphism between $\uht * L$ and Khovanov's original $sl_3$-link homology, 
which in our notation is equal to $\uhtabc 0 0 0 * L$.

In the following two subsections we study the cases in which $f(X)$ has two or 
three distinct roots. We first work out the case for three distinct roots, 
because this case has been done already by Gornik~\cite{gornik} essentially. 
Even in this case we define and prove 
everything precisely and completely. We have two good reasons for doing this. 
First of all 
we generalize Gornik's work to the arbitrary case of three distinct roots, 
whereas he, strictly speaking, only considers the particular case of the 
third roots of unity. Given the definitions and arguments for the general 
case, one easily recognizes Gornik's definitions and arguments for his 
particular case. Working one's way back is harder, 
also because Gornik followed the approach using matrix factorizations and not 
cobordisms. Secondly these general definitions and arguments are necessary 
for understanding the last subsection, where we treat the case in which 
$f(X)$ has only two distinct roots, which is clearly different from Gornik's.


\subsection{Three distinct roots}\label{sec:3roots}

In this subsection we assume that 
the three roots of $f(X)$, denoted $\alpha,\beta,\gamma\in\bC$, are all 
distinct.
First we determine Gornik's idempotents in the algebra 
$\bC[X]/\left(f(X)\right)$. By the Chinese Remainder Theorem we have the 
following isomorphism of 
algebras
$$\bC[X]/\left(f(X)\right)\cong \bC[X]/\left(X-\alpha\right)\oplus 
\bC[X]/\left(X-\beta\right)\oplus \bC[X]/\left(X-\gamma\right)\cong \bC^3.$$

\begin{defn} Let $Q_{\alpha}(X)$, $Q_{\beta}(X)$ e $Q_{\gamma}(X)$ be the 
idempotents 
in 
$\bC[X]/\left(f(X)\right)$ corresponding to $(1,0,0)$, $(0,1,0)$ e $(0,0,1)$ 
in $\bC^3$ under the isomorphism in the Chinese Remainder Theorem.
\end{defn} 
\noindent As a matter of fact it is easy to compute the idempotents 
explicitly: 
$$Q_{\alpha}(X)=\dfrac{(X-\beta)(X-\gamma)}{(\alpha-\beta)(\alpha-\gamma)},\, 
Q_{\beta}(X)=\dfrac{(X-\alpha)(X-\gamma)}{(\beta-\alpha)(\beta-\gamma)},$$
$$Q_{\gamma}(X)=\dfrac{(X-\alpha)(X-\beta)}{(\gamma-\alpha)(\gamma-\beta)}.$$

\noindent By definition we get 
\begin{lem}
\label{lem:chinese}
$$Q_{\alpha}(X)+Q_{\beta}(X)+Q_{\gamma}(X)=1,$$
$$Q_{\alpha}(X)Q_{\beta}(X)=Q_{\alpha}(X)Q_{\gamma}(X)=Q_{\beta}(X)
Q_{\gamma}(X)=0,$$
$$Q_{\alpha}(X)^2=Q_{\alpha}(X),\quad Q_{\beta}(X)^2=Q_{\beta}(X),\quad 
Q_{\gamma}(X)^2=Q_{\gamma}(X).$$
\end{lem}

Let $\Gamma$ be a resolution of a link L and let $E(\Gamma)$ be the set of all 
edges in $\Gamma$. 
In \cite{khovanovsl3} Khovanov defines the following algebra (in his case 
for $a=b=c=0$). 
\begin{defn} Let $R(\Gamma)$ be the free commutative algebra generated by 
the elements $X_i$, with $i\in E(\Gamma)$, modulo the relations 
\begin{equation}
\label{eq:relations}
X_i+X_j+X_k=a,\quad X_iX_j+X_jX_k+X_iX_k=-b,\quad X_iX_jX_k=c,
\end{equation}
for any triple of edges $i,j,k$ which share a trivalent vertex.
\end{defn}

The following definitions and results are analogous to Gornik's results in 
sections 2 and 3 of \cite{gornik}.
Let $S=\left\{\alpha,\beta,\gamma\right\}$.  

\begin{defn}
\label{defn:coloring}
A {\em coloring} of $\Gamma$ is defined to be a map 
$\phi\colon E(\Gamma)\to S$. 
Denote the set of all colorings by $S(\Gamma)$.
An {\em admissible coloring} is a coloring such that 
\begin{equation}
\label{eq:admissible}
\begin{array}{lll}
a&=&\phi(i)+\phi(j)+\phi(k)\\
-b&=&\phi(i)\phi(j)+\phi(j)\phi(k)+
\phi(i)\phi(k)\\
c&=&\phi(i)\phi(j)\phi(k),
\end{array} 
\end{equation}
for any edges $i,j,k$ incident to the same trivalent vertex. 
Denote the set of all admissible colorings by 
$AS(\Gamma)$. 
\end{defn}
\noindent Of course admissibility 
is equivalent to requiring that the three colors $\phi(i),\phi(j)$ and 
$\phi(k)$ be all distinct. 

A simple calculation shows that $f(X_i)=0$ in $R(\Gamma)$, for any 
$i\in E(\Gamma)$. Therefore, for any edge $i\in E(\Gamma)$, there exists a 
homomorphism of algebras from 
$\bC[X]/(F(X))$ to 
$R(\Gamma)$ defined by $X\mapsto X_i$. 
Thus, for any coloring $\phi$, we define
\begin{defn} 
\label{defn:idempotents}
$$Q_{\phi}(\Gamma)=\prod_{i\in E(\Gamma)}Q_{\phi(i)}(X_i)\in R(\Gamma).$$
\end{defn}
\noindent Lemma~\ref{lem:chinese} implies the following corollary. 
\begin{cor}
\label{cor:chinese}
$$\sum_{\phi\in S(\Gamma)}Q_{\phi}(\Gamma)=1,$$
$$Q_{\phi}(\Gamma)Q_{\psi}(\Gamma)=\delta_{\psi}^{\phi}Q_{\phi},$$
where $\delta_{\psi}^{\phi}$ is the Kronecker delta. 
\end{cor}

\noindent Note that the definition of $Q_{\phi}(\Gamma)$ implies that 
\begin{equation}
\label{eq:eigenvalue}
X_iQ_{\phi}(\Gamma)=\phi(i)Q_{\phi}(\Gamma).
\end{equation}

The following lemma is our analogue of Gornik's theorem 3. 

\begin{lem} 
\label{lem:decomp}
For any non-admissible coloring $\phi$, we have
$$Q_{\phi}(\Gamma)=0.$$ 

For any admissible coloring $\phi$, we have
$$Q_{\phi}(\Gamma)R(\Gamma)\cong \bC.$$

Therefore, we get a direct sum decomposition 
$$R(\Gamma)\cong \bigoplus_{\phi\in AS(\Gamma)}\bC Q_{\phi}(\Gamma).$$
\end{lem}
\begin{proof} Let $\phi$ be any coloring and let $i,j,k\in E(\Gamma)$ be 
three edges sharing a trivalent vertex. By the relations in  
(\ref{eq:relations}) and equation (\ref{eq:eigenvalue}), we get 
\begin{equation}
\label{eq:releigenvalues}
\begin{array}{lll}
aQ_{\phi}(\Gamma)&=&(\phi(i)+\phi(j)+\phi(k))Q_{\phi}(\Gamma)\\
-bQ_{\phi}(\Gamma)&=&(\phi(i)\phi(j)+\phi(j)\phi(k)+
\phi(i)\phi(k))Q_{\phi}(\Gamma)\\
cQ_{\phi}(\Gamma)&=&\phi(i)\phi(j)\phi(k)Q_{\phi}(\Gamma).
\end{array} 
\end{equation}
If $\phi$ is non-admissible, then, by comparing (\ref{eq:admissible}) and 
(\ref{eq:releigenvalues}), we see that $Q_{\phi}(\Gamma)$ vanishes.

Now suppose $\phi$ is admissible. Recall that $R_{\phi}(\Gamma)$ is a 
quotient of the algebra 
\begin{equation}
\label{eq:algebra}
\bigotimes_{i\in E(\Gamma)}\bC[X_i]/
\left(f(X_i)\right).
\end{equation} 
Just as in definition~\ref{defn:idempotents} we can define the idempotents 
in the algebra in (\ref{eq:algebra}), which we also denote 
$Q_{\phi}(\Gamma)$. By the 
Chinese Remainder Theorem, there is a projection of the algebra in 
(\ref{eq:algebra}) onto 
$\bC$, which maps 
$Q_{\phi}(\Gamma)$ to $1$ and $Q_{\psi}(\Gamma)$ to $0$, for any 
$\psi\ne\phi$. It is not hard to see that, since $\phi$ is admissible, 
that projection factors through the quotient $R(\Gamma)$, which implies the 
second claim in the lemma. 
\end{proof}

\noindent As in \cite{khovanovsl3}, the relations in figure~\ref{fig:pdots} show that 
$R(\Gamma)$ acts on $C(\Gamma)$ by the usual action induced by the cobordism 
which merges a circle and the relevant edge of $\Gamma$. Let us write 
$C_{\phi}(\Gamma)=Q_{\phi}(\Gamma)C(\Gamma)$. By 
corollary~\ref{cor:chinese} and lemma~\ref{lem:decomp}, we have  
a direct sum decomposition 
\begin{equation}
\label{eq:decomp}
C(\Gamma)=\bigoplus_{\phi\in AS(\Gamma)}C_{\phi}(\Gamma).
\end{equation}
Note that we have
\begin{equation}
\label{eq:condeigenvalues}
z\in C_{\phi}(\Gamma)\Leftrightarrow \forall i,\, X_iz=\phi(i)z
\end{equation}
for any $\phi\in AS(\Gamma)$.

Let $\phi$ be a coloring of the arcs of $L$ by $\alpha,\beta$ and $\gamma$. 
Note that $\phi$ induces a unique coloring of the unmarked edges of any 
resolution of $L$. 
\begin{defn}
\label{defn:canonical}
We say that a coloring of the arcs of $L$ is {\em admissible} if there 
exists a resolution of $L$ which admits a compatible admissible coloring. 
Note that if such a resolution exists, its coloring is uniquely determined 
by $\phi$, so we use the same notation. Note also that an admissible 
coloring of $\phi$ induces a unique admissible coloring of $L$. If $\phi$ 
is an admissible coloring, we call the elements in $C_{\phi}(\Gamma)$ 
{\em admissible cochains}. We denote the set of all admissible colorings 
of $L$ by $AS(L)$.

We say that an admissible coloring of 
$L$ is a {\em canonical coloring} if the arcs belonging to the same component 
of $L$ have the same color. If $\phi$ is a canonical coloring, we call 
the elements in $C_{\phi}(\Gamma)$ {\em canonical cochains}. 
We denote the set of canonical colorings of $L$ by $S(L)$.   
\end{defn}

Note that, for a fixed $\phi\in AS(L)$, the admissible 
cochain groups $C_{\phi}(\Gamma)$ form a subcomplex 
$\cut * {L_{\phi}}\subseteq \cut * {L}$ whose homology we denote by 
$\uht * {L_{\phi}}$. The following lemma shows that only the canonical 
cochain groups matter, as Gornik indicated in his remarks before his 
main theorem 2 in~\cite{gornik}. 

\begin{thm} 
\label{thm:decomp}
$$\uht * L=\bigoplus_{\phi\in S(L)}\uht * {L_{\phi}}.$$
\end{thm}
\begin{proof}
By (\ref{eq:decomp}) we have 
$$\uht * L=\bigoplus_{\phi\in AS(L)}\uht * {L_{\phi}}.$$
\begin{figure}[h]
$$\xymatrix@R=1.8mm{
\includegraphics[height=0.7in]{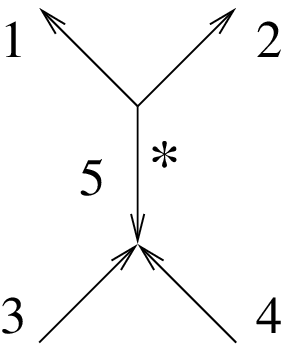} &
\includegraphics[height=0.7in]{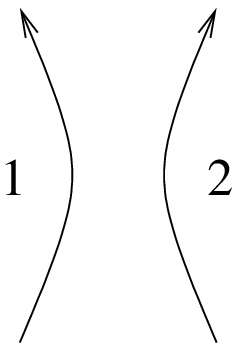} \\
\Gamma & \Gamma'
}$$
\caption{Ordering edges}
\label{fig:edge_order_tvert}
\end{figure}
Let us now show that $\uht * {L_{\phi}}=0$ if $\phi$ is admissible but 
non-canonical. 
Let $\Gamma$ and $\Gamma'$ be the diagrams in figure~\ref{fig:edge_order_tvert}, which are the 
boundary of the cobordism which defines the differential in $\cut * L$, 
and order their edges as indicated. Up to permutation, the 
only admissible colorings of $\Gamma$ 
are
$$\xymatrix@R=1.8mm{
\includegraphics[height=0.7in]{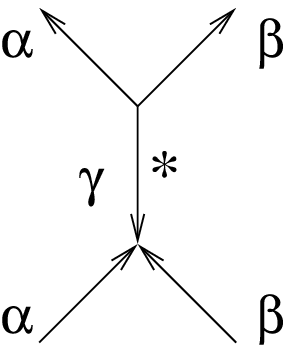} &\text{and}&
\includegraphics[height=0.7in]{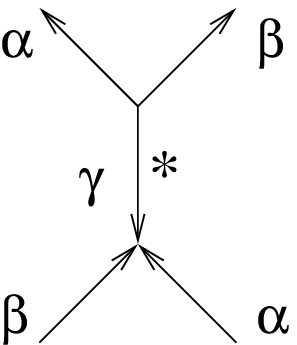} \\
\phi_1 & & \phi_2
}.$$
Up to permutation, the only admissible colorings of $\Gamma'$ 
are 
$$\xymatrix@R=1.8mm{
\includegraphics[height=0.7in]{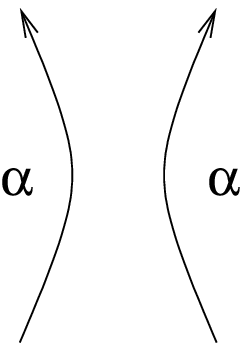} &\text{and}&
\includegraphics[height=0.7in]{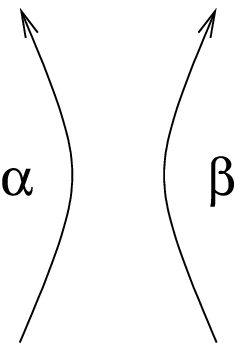} \\
\phi_0' & & \phi_1' 
}.$$
Note that only $\phi_2$ and $\phi_0'$ can be canonical. 
We get  
\begin{equation}
\label{eq:canonical}
\begin{array}{lll}
0&\leftarrow&C_{\phi'_0}(\Gamma'),\\
C_{\phi_1}(\Gamma)&\cong& C_{\phi_1'}(\Gamma'),\\
C_{\phi_2}(\Gamma)&\rightarrow& 0.\\
\end{array}
\end{equation}
Note that the elementary cobordism has to map colorings to compatible colorings. This 
explains the first and the third line. Let us explain the second line. Apply 
the elementary cobordisms $\Gamma'\to \Gamma\to \Gamma'$ and 
use relation (\emph{RD}) on 
page~\pageref{R-rel} to obtain the linear map 
$C_{\phi_{1,2}'}(\Gamma')\to C_{\phi_{1,2}'}(\Gamma')$ 
given by 
$$z\mapsto (\beta-\alpha)z.$$ 
Since $\alpha\ne\beta$, we see that this map is injective. Therefore 
the map $C_{\phi_{1,2}'}(\Gamma')\to C_{\phi_{1,2}}(\Gamma)$ is injective 
too. A similar argument, using the (\emph{DR}) relation, shows that 
$C_{\phi_{1,2}}(\Gamma)\to C_{\phi_{1,2}'}(\Gamma')$ is injective. Therefore 
both maps are isomorphisms. 

Next, let $\phi$ be admissible but non-canonical. 
Then there exists at least one 
crossing, denoted $c$, in $L$ which has a resolution with a non-canonical coloring. 
Let 
$\cut * {L_{\phi}^1}$ be the subcomplex of $\cut * {L_{\phi}}$ defined by the 
resolutions of $L$ in which $c$ has been resolved by the $1$-resolution. 
Let $\cut * {L_{\phi}^0}$ be the complex obtained 
from $\cut * {L_{\phi}}$ by deleting all resolutions which do not belong to 
$\cut * {L_{\phi}^1}$ and all arrows which have a source or target which is 
not one of the remaining resolutions. Note that we have a short exact sequence 
of complexes 
\begin{equation}
\label{eq:ses}
0\to\cut * {L_{\phi}^1}\to\cut * {L_{\phi}}\to\cut * {L_{\phi}^0}\to 0.
\end{equation}
The isomorphism in (\ref{eq:canonical}) shows that the natural map 
$$ \cut * {L_{\phi}^0}\to \cut {* + 1} {L_{\phi}^1},$$
defined by the elementary cobordisms which induce the 
connecting homomorphism in the long exact sequence associated to 
(\ref{eq:ses}), is an isomorphism. 
By exactness of this long exact sequence we see that 
$\uht * {L_{\phi}}=0$.
\end{proof}

\begin{lem}
\label{lem:cyclic}
For any $\phi\in AS(\Gamma)$, we have 
$C_{\phi}(\Gamma)\cong\bC$.
\end{lem}
\begin{proof} 
We use induction with respect to $v$, the number of trivalent vertices in 
$\Gamma$.
The claim is obviously true for a circle. Suppose $\Gamma$ has a digon, 
with the edges ordered as in figure~\ref{fig:edge_order_digon}. 
\begin{figure}[h]
$$\xymatrix@R=1.8mm{
\includegraphics[height=0.7in]{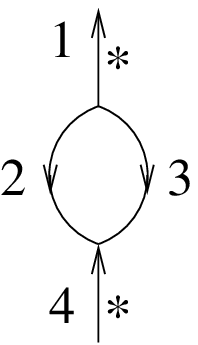} &
\includegraphics[height=0.7in]{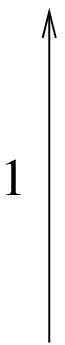} \\
\Gamma & \Gamma'
}$$
\caption{Ordering edges in digon}
\label{fig:edge_order_digon}
\end{figure}
Note 
that $X_1=X_4\in R(\Gamma)$ holds as a consequence of the relations 
in (\ref{eq:relations}). Let 
$\Gamma'$ be the web obtained by removing the digon, as in figure~\ref{fig:edge_order_digon}. 
Up to permutation, the only possible admissible colorings of $\Gamma$ and the 
corresponding admissible coloring of $\Gamma'$ are 
$$\xymatrix@R=1.8mm{
\includegraphics[height=0.7in]{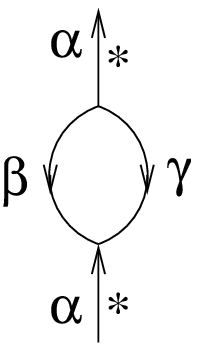} &
\includegraphics[height=0.7in]{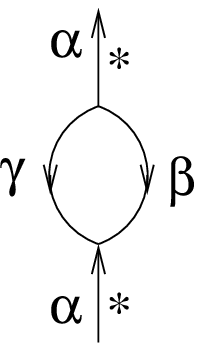} &
\includegraphics[height=0.7in]{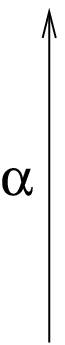}
\\
\phi_1 & \phi_2 & \phi' 
}.$$
The \emph{Digon Removal} isomorphism in 
lemma~\ref{lem:KhK} yields
$$C_{\phi_1}(\Gamma)\oplus C_{\phi_2}(\Gamma)
\cong C_{\phi'}(\Gamma')\oplus 
C_{\phi'}(\Gamma').$$
By induction, we have $C_{\phi'}(\Gamma')\cong \bC$, so 
$\dim C_{\phi_1}(\Gamma)+\dim C_{\phi_2}(\Gamma)=2$. 
For symmetry reasons this implies that 
$\dim C_{\phi_1}(\Gamma)=\dim C_{\phi_2}(\Gamma)=1$, 
which proves the claim. To be a bit more precise, let 
$B_{\beta,\gamma}$ and $B_{\gamma,\beta}$ be the following two colored 
cobordisms:
$$B_{\beta,\gamma}=
\raisebox{-16pt}{\includegraphics[height=0.66in]{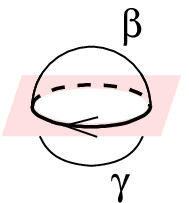}}
, \qquad 
B_{\gamma,\beta}=
\raisebox{-16pt}{\includegraphics[height=0.6in]{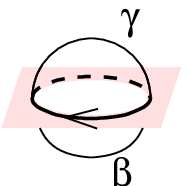}}.$$ 
Note that we have 
$$B_{\beta,\gamma}+B_{\gamma,\beta}=0$$ by
$$\raisebox{-12pt}{\includegraphics[height=0.4in]{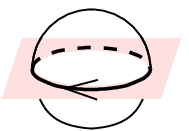}}=0,$$ 
and by
$$\raisebox{-12pt}{\includegraphics[height=0.4in]{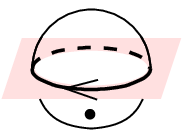}}=
\raisebox{-5pt}{\includegraphics[height=0.21in,width=0.52in]{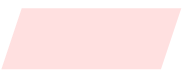}}
$$ we have
$$\gamma B_{\beta,\gamma}+\beta B_{\gamma,\beta}=
\mbox{id}_{C_{\phi'}(\Gamma')}.$$ 
These two identities imply 
$$(\gamma-\beta)B_{\beta,\gamma}=(\beta-\gamma)B_{\gamma,\beta}=
\mbox{id}_{C_{\phi'}(\Gamma')}.$$ 
Therefore we conclude
that $C_{\phi_1}(\Gamma)$ and $C_{\phi_2}(\Gamma)$ are non-zero, 
which for dimensional reasons implies 
$\dim C_{\phi_1}(\Gamma)=\dim C_{\phi_2}(\Gamma)=1$.

Now, suppose $\Gamma$ contains a square, with the edges ordered as in 
figure~\ref{fig:edge_order_square} left. 
\begin{figure}[h]
$$\xymatrix@R=1.8mm{
 \includegraphics[height=0.7in]{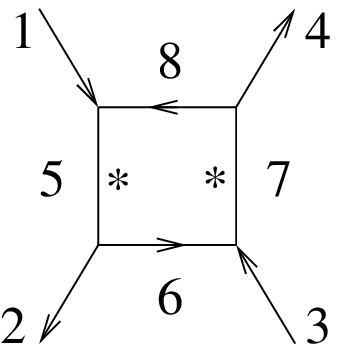} &
 \includegraphics[height=0.7in]{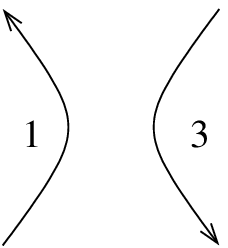} &
 \includegraphics[height=0.7in]{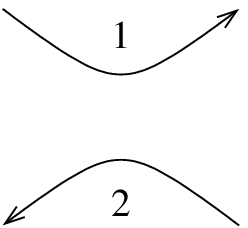} \\
 \Gamma & \Gamma' & \Gamma''
}$$
\caption{Ordering edges in square}
\label{fig:edge_order_square}
\end{figure}
Let $\Gamma'$ and $\Gamma''$ be the two corresponding webs under the \emph{Square Removal} isomorphism in lemma~\ref{lem:KhK}. Up to permutation there is only one admissible 
non-canonical coloring and one canonical coloring:
$$
\xymatrix@R=1.8mm
{\includegraphics[height=0.7in]{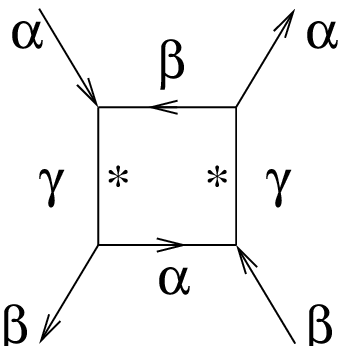} & 
\includegraphics[height=0.7in]{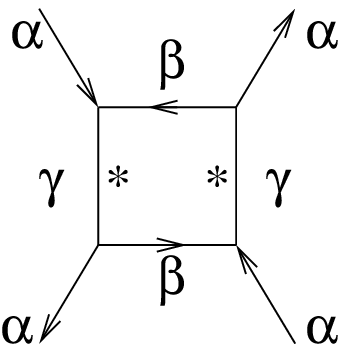} \\
\text{canonical} & \text{non-canonical}
}.$$

Let us first consider the canonical coloring. Clearly $C_{\phi}(\Gamma)$ 
is isomorphic to $C_{\phi''}(\Gamma'')$, where $\phi''$ is the 
unique compatible canonical coloring, because there is no 
compatible coloring of $\Gamma'$. Therefore the result follows by induction.

Now consider the admissible non-canonical coloring. As proved in 
theorem~\ref{thm:decomp} we have the following isomorphism:
$$\raisebox{-22pt}{\includegraphics[height=0.7in]{sqcol_aaaagbgb}}
\cong
\raisebox{-22pt}{\includegraphics[height=0.7in]{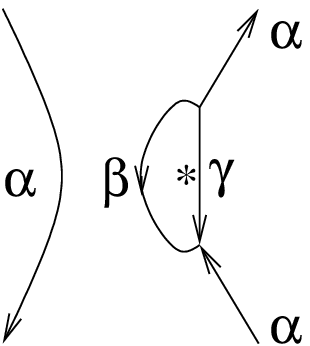}}.$$
By induction the right-hand side is one-dimensional, which proves the claim. 

\end{proof}

Thus we arrive at Gornik's main theorem 2. Note that there are $3^n$ 
canonical colorings of $L$, where $n$ is the number of components of $L$. 
Note also that the homological degrees
of the canonical cocycles are easy to compute, because we know that the 
canonical cocycles corresponding to the oriented resolution without vertices 
have homological degree zero. 

\begin{thm}
\label{thm:gornik} 
The dimension of $\uht * L$ equals $3^n$, where $n$ is the number of 
components of $L$.  

For any $\phi\in S(L)$, there exists a non-zero element 
$a_{\phi}\in \uht i L$, unique up to a scalar, where 
$$i=\sum_{(\epsilon_1,\epsilon_2)\in S\times S,\, \epsilon_1\ne \epsilon_2}
\lk(\phi^{-1}(\epsilon_1),\phi^{-1}(\epsilon_2)).$$
\end{thm}


\subsection{Two distinct roots}\label{sec:2roots}

In this section we assume that $f(X)=(X-\alpha)^2(X-\beta)$, with $\alpha\ne 
\beta$.
We follow an approach similar to the one in the previous 
section. First we define the relevant idempotents. By the Chinese Remainder 
Theorem we have 
$$\bC[X]/(f(X))\cong \bC[X]/((X-\alpha)^2)\oplus \bC[X]/(X-\beta).$$

\begin{defn} Let $Q_{\alpha}$ and $Q_{\beta}$  be the idempotents 
in $\bC[X]/(f(X))$ corresponding to $(1,0)$ and $(0,1)$ in  
$\bC[X]/((X-\alpha)^2)\oplus \bC[X]/(X-\beta)$ under the above isomorphism.   
\end{defn}

\noindent Again it is easy to compute the idempotents explicitly:
$$Q_{\alpha}=1-\dfrac{(X-\alpha)^2}{(\beta-\alpha)^2},
\quad Q_{\beta}=\dfrac{(X-\alpha)^2}
{(\beta-\alpha)^2}.$$

\noindent By definition we get 

\begin{lem}
\label{lem:chinese2}
$$Q_{\alpha}+Q_{\beta}=1,$$
$$Q_{\alpha}Q_{\beta}=0,$$
$$Q_{\alpha}^2=Q_{\alpha},\quad Q_{\beta}^2=Q_{\beta}.$$
\end{lem}

\noindent Throughout this subsection let $S=\{\alpha,\beta\}$. 
We define colorings of webs and admissibility as in 
definition~\ref{defn:coloring}. Note that a coloring 
is admissible if and only if at each trivalent vertex the, unordered, incident 
edges are colored $\alpha,\alpha,\beta$. Let $\Gamma$ be a web and $\phi$ a 
coloring. The definition of the 
idempotents $Q_{\phi}(\Gamma)$ in $R(\Gamma)$ is the same as in 
definition~\ref{defn:idempotents}. Clearly, corollary~\ref{cor:chinese} also 
holds in this section. 
However, equation~\ref{eq:eigenvalue} changes. By the Chinese Remainder 
Theorem, we get  
\begin{equation}
\label{eq:eigenvalue2}
\begin{cases}
(X_i-\beta)Q_{\phi}(\Gamma)=0,&\text{if}\quad \phi(i)=\beta\\
(X_i-\alpha)^2Q_{\phi}(\Gamma)=0,&\text{if}\quad \phi(i)=\alpha.
\end{cases}
\end{equation}

Lemma~\ref{lem:decomp} also changes. Its analogue becomes:

\begin{lem}
\label{lem:decomp2}
For any non-admissible coloring $\phi$, we have
$$Q_{\phi}(\Gamma)=0.$$ 
Therefore, we have a direct sum decomposition 
$$R(\Gamma)\cong \bigoplus_{\phi\in AS(\Gamma)} Q_{\phi}(\Gamma)R(\Gamma).$$

For any $\phi\in AS(\Gamma)$, we have  
$\dim Q_{\phi}(\Gamma)R(\Gamma)=2^m$, where $m$ 
is the number of cycles in $\phi^{-1}(\alpha)\subseteq 
\Gamma$.
\end{lem}
\begin{proof}
First we prove that inadmissible colorings yield trivial idempotents. 
Let $\phi$ be any coloring of 
$\Gamma$ and let $i,j,k$ be three edges sharing a trivalent vertex. First 
suppose that all edges are colored by $\beta$. By equations 
(\ref{eq:eigenvalue2}) we 
get 
$$aQ_{\phi}(\Gamma)=(X_i+X_j+X_k)Q_{\phi}(\Gamma)=3\beta Q_{\phi}(\Gamma),$$
which implies that $Q_{\phi}(\Gamma)=0$, because $a=2\alpha+\beta$ and 
$\alpha\ne\beta$.

Next suppose $\phi(i)=\phi(j)=\beta$ and $\phi(k)=\alpha$. Then 
$$aQ_{\phi}(\Gamma)=(X_i+X_j+X_k)Q_{\phi}(\Gamma)=(2\beta+X_k) 
Q_{\phi}(\Gamma).$$
Thus $X_kQ_{\phi}(\Gamma)=(2\alpha-\beta)Q_{\phi}(\Gamma)$. Therefore we get  
$$0=(X_k-\alpha)^2Q_{\phi}(\Gamma)=(\alpha-\beta)^2Q_{\phi}(\Gamma),$$
which again implies that $Q_{\phi}(\Gamma)=0$. 

Finally, suppose $i,j,k$ are all colored by $\alpha$. Then we have 
$$\left((X_i-\alpha)^2+(X_j-\alpha)^2+(X_k-\alpha)^2\right)
Q_{\phi}(\Gamma)=0.$$
Using the relations in (\ref{eq:relations}) we get
$$(X_i-\alpha)^2+(X_j-\alpha)^2+(X_k-\alpha)^2=(\alpha-\beta)^2,$$
so we see that $Q_{\phi}(\Gamma)=0$.

Now, let $\phi$ be an admissible coloring. Note that the admissibility 
condition implies that $\phi^{-1}(\alpha)$ consists of a disjoint union of 
cycles. To avoid confusion, let us remark that we do not take into consideration 
the orientation of the edges when we speak about cycles, as one would in algebraic 
topology. What we mean by a cycle is simply a piece-wise linear closed loop. 
Recall that $R(\Gamma)$ is a quotient of the algebra 
\begin{equation}
\label{eq:algebra2}
\bigotimes_{i\in E(\Gamma)}\bC[X_i]/\left(f(X_i)\right)
\end{equation}
and that we can define idempotents, also denoted 
$Q_{\phi}(\Gamma)$, in the latter. 
Note that by the Chinese Remainder Theorem 
there exists a homomorphism of algebras which projects the algebra in 
(\ref{eq:algebra2})  
onto 
\begin{equation}
\label{eq:algebra3}
\bigotimes_{\phi(i)=\alpha}\bC[X_i]/\left((X_i-\alpha)^2\right)\bigotimes 
\bigotimes_{\phi(i)=\beta}\bC[X_i]/\left(X_i-\beta\right),
\end{equation}
which maps $Q_{\phi}(\Gamma)$ to $1$ and $Q_{\psi}(\Gamma)$ to $0$, for any 
$\psi\ne \phi$. 
Define 
$R_{\phi}(\Gamma)$ to be the quotient of the algebra in (\ref{eq:algebra3}) 
by the relations $X_i+X_j=2\alpha$, for all edges $i$ and $j$ 
which share a trivalent vertex and satisfy $\phi(i)=\phi(j)=\alpha$. 
Note that $X_iX_j=\alpha^2$ also holds in $R_{\phi}(\Gamma)$, for such edges $i$ 
and $j$. 

Suppose that the edges $i,j,k$ are 
incident to a trivalent vertex in $\Gamma$ and that they are colored 
$\alpha,\alpha,\beta$. It is easy to see that by the projection onto 
$R_{\phi}(\Gamma)$ we get 
$$
\begin{array}{lll}
X_i+X_j+X_k&\mapsto& a\\
X_iX_j+X_iX_k+X_jX_k&\mapsto& -b\\
X_iX_jX_k&\mapsto&c.
\end{array}
$$ 
Therefore the projection descends to a projection from $R(\Gamma)$ onto 
$R_{\phi}(\Gamma)$. Since $Q_{\phi}(\Gamma)$ is mapped to $1$ and 
$Q_{\psi}(\Gamma)$ to $0$, 
for all $\psi\ne \phi$, we see that the projection restricts to a 
surjection of algebras  
$$Q_{\phi}(\Gamma)R(\Gamma)\to R_{\phi}(\Gamma).$$ A simple computation 
shows that the equality  
$$X_i+X_jQ_{\phi}(\Gamma)=2\alpha Q_{\phi}(\Gamma)$$
holds in $R(\Gamma)$, which implies that the surjection above is an isomorphism 
of algebras. This proves the final claim in the lemma.  
\end{proof}

As in (\ref{eq:condeigenvalues}), 
for any $\phi\in AS(\Gamma)$, we get 
\begin{equation}
\label{eq:condeigenvalues2}
z\in C_{\phi}(\Gamma)\Leftrightarrow 
\begin{cases}(X_i-\beta)z=0,&\forall i\colon \phi(i)=\beta,\\
(X_i-\alpha)^2z=0,&\forall i\colon \phi(i)=\alpha.
\end{cases}
\end{equation}

\noindent Let $\phi$ be a coloring of the arcs of $L$ by $\alpha$ and $\beta$. 
Note that $\phi$ induces a unique coloring of the unmarked edges of any 
resolution of $L$. We define {\em admissible} and {\em canonical colorings} 
of $L$ as in definition~\ref{defn:canonical}.

Note, as before, that, for a fixed admissible coloring $\phi$ of $L$, the 
admissible 
cochain groups $C_{\phi}(\Gamma)$ form a subcomplex 
$\cut * {L_{\phi}}\subseteq \cut * {L}$ whose homology we denote by 
$\uht * {L_{\phi}}$. The following theorem is the analogue of 
theorem~\ref{thm:decomp}. 

\begin{thm} 
\label{thm:decomp2}
$$\uht * L=\bigoplus_{\phi\in S(L)}\uht * {L_{\phi}}.$$
\end{thm}
\begin{proof}
By lemma~\ref{lem:decomp2} we get 
$$\uht * L=\bigoplus_{\phi\in AS(L)}\uht * {L_{\phi}}.$$

Let us now show that $\uht * {L_{\phi}}=0$ if $\phi$ is admissible but 
non-canonical. 
Let $\Gamma$ and $\Gamma'$ be the diagrams in figure~\ref{fig:edge_order_tvert}, which are the 
boundary of the cobordism which induces the differential in $\cut * L$, 
and order their edges as indicated. The only admissible colorings of $\Gamma$ 
are 
$$\xymatrix@R=1.8mm{
\includegraphics[height=0.7in]{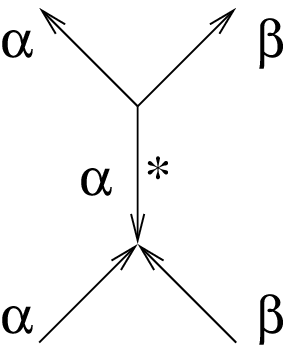} &
\includegraphics[height=0.7in]{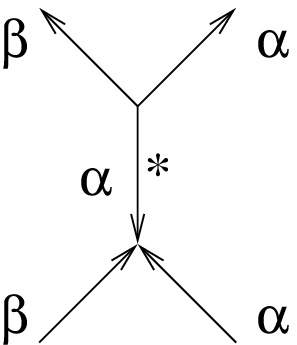} &
\includegraphics[height=0.7in]{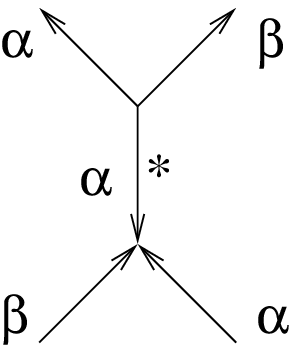} &
\includegraphics[height=0.7in]{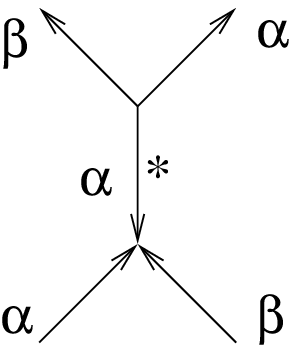} &
\includegraphics[height=0.7in]{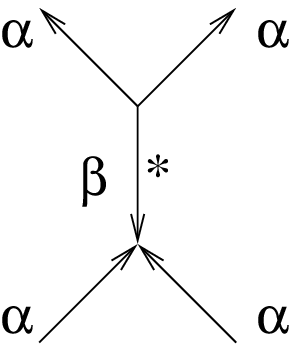}
\\
\phi_1 & \phi_2 & \phi_3 & \phi_4 & \phi_5 
}.$$
The only admissible colorings of $\Gamma'$ 
are 
$$\xymatrix@R=1.8mm{
\includegraphics[height=0.7in]{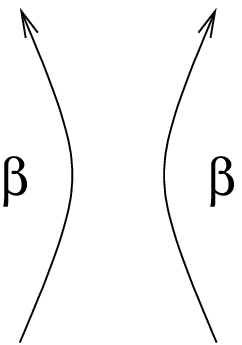} &
\includegraphics[height=0.7in]{tedgescol_ab} &
\includegraphics[height=0.7in]{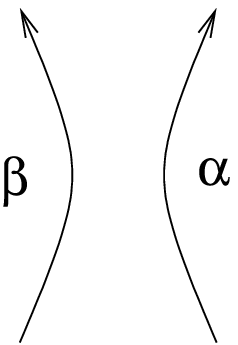} &
\includegraphics[height=0.7in]{tedgescol_aa}   
\\
\phi_0' & \phi_1' & \phi_2' & \phi_5' 
}.$$
Note that only $\phi_3,\phi_4, \phi_5, \phi_0'$ and $\phi_5'$ can be 
canonical. We get
\begin{equation}
\label{eq:canonical2}
\begin{array}{lll}
0&\leftarrow&C_{\phi'_0}(\Gamma'),\\
C_{\phi_1}(\Gamma)&\cong& C_{\phi_1'}(\Gamma'),\\
C_{\phi_2}(\Gamma)&\cong& C_{\phi_2'}(\Gamma'),\\
C_{\phi_3}(\Gamma)&\rightarrow& 0,\\
C_{\phi_4}(\Gamma)&\rightarrow& 0,\\
C_{\phi_5}(\Gamma)&\leftrightarrow& C_{\phi_5'}(\Gamma').
\end{array}
\end{equation}
Note that the last line in the list above only states that the cobordism  
induces a map from one side to the other or vice-versa, 
but not that it is an isomorphism in 
general. The second and third line contain isomorphisms. Let us explain the 
second line, the third being similar. Apply 
the elementary cobordism $\Gamma'\to\Gamma\to\Gamma'$ and 
use relation (\emph{RD}) on 
page~\pageref{R-rel} to obtain the linear map 
$C_{\phi_{1,2}'}(\Gamma')\to C_{\phi_{1,2}'}(\Gamma')$ 
given by 
$$z\mapsto (\beta-X_1)z.$$ 
Suppose $(X_1-\beta)z=0$. Then $z\in C_{\phi_{0}'}(\Gamma')$, 
because 
$X_1z=X_2z=\beta z$. This implies that $z\in C_{\phi_{1,2}'}(\Gamma')
\cap C_{\phi_{0}'}(\Gamma')=\{0\}$. Thus the 
map above is injective, and therefore the map 
$C_{\phi_{1,2}'}(\Gamma')\to C_{\phi_{1,2}}(\Gamma)$ is injective. 
A similar argument, using the relation (\emph{DR}), shows that 
$C_{\phi_{1,2}}(\Gamma)\to C_{\phi_{1,2}'}(\Gamma')$ is injective. Therefore 
both maps are isomorphisms. 

The isomorphisms in (\ref{eq:canonical2}) imply that $\uht * {L_{\phi}}=0$ holds, 
when $\phi$ is admissible but non-canonical, as we explained in the proof 
of theorem~\ref{thm:decomp}.
\end{proof}

Let $C_{\phi}(\Gamma)$ be a canonical cochain group. In this case it does not 
suffice to compute the dimensions of $C_{\phi}(\Gamma)$, for all $\phi$ and 
$\Gamma$, because we also need to determine the differentials. Therefore we first 
define a canonical cobordism in $C_{\phi}(\Gamma)$. 

\begin{figure}[h]
\includegraphics[height=0.2in]{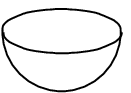} \qquad
\includegraphics[height=0.4in]{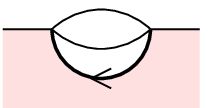}\qquad
\includegraphics[height=0.7in]{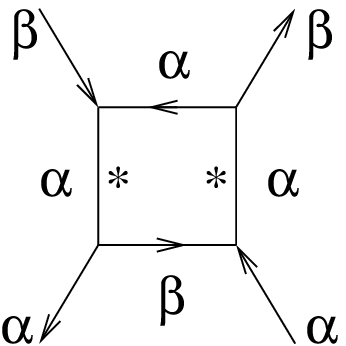}
\raisebox{0.3in}{:}\quad
\vspace{2ex}
\includegraphics[height=0.7in]{ssquareup}
$$
\raisebox{-0.3in}{\includegraphics[height=0.7in]{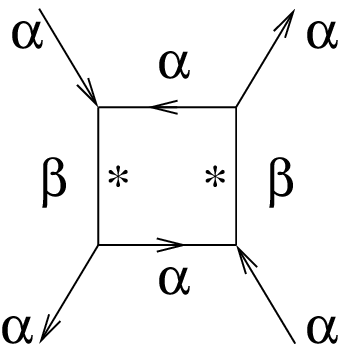}}:
\left\{
\begin{matrix}
 \includegraphics[height=0.7in]{ssquareup} & 
 \raisebox{0.4in}{\text{if }\;\raisebox{-0.2pt}{\includegraphics[height=3.2mm]{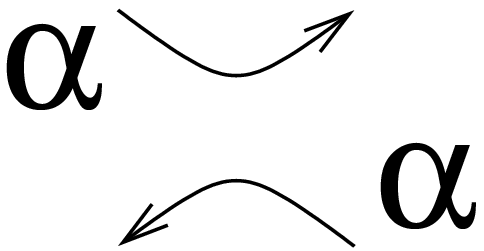}}\;\text{ belong to different $\alpha$-cycles}}  \\
\includegraphics[height=0.7in]{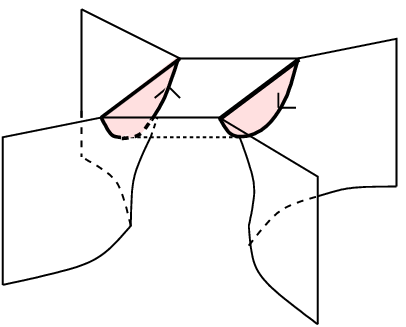} & 
 \raisebox{0.4in}{\text{if }\;\raisebox{-0.2pt}{\includegraphics[height=3.2mm]{tedgesop}}\;\text{ belong to same $\alpha$-cycle}}
\end{matrix}\right.
$$
\caption{Elementary cobordisms}
\label{fig:ecob_x3=x2}
\end{figure}

\begin{figure}[h]
\includegraphics[height=0.35in]{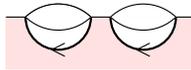}
\caption{several digons}
\label{fig:sev-dig}
\end{figure}

\begin{defn}
Let $\phi\in S(L)$. We define a cobordism 
$\Sigma_{\phi}(\Gamma)\colon \emptyset\to\Gamma$ by 
gluing together the elementary cobordisms in figure~\ref{fig:ecob_x3=x2} and 
multiplying by $Q_{\phi}(\Gamma)$. 
We call $\Sigma_{\phi}(\Gamma)$ the {\em canonical cobordism} 
in $C_{\phi}(\Gamma)$. 
\end{defn}
\noindent For any canonically colored web, we can find a way to build up the canonical 
cobordism using only the above elementary cobordisms with canonical colorings, except 
when we have several digons as in figure~\ref{fig:sev-dig} where we might 
have to stick in two digons at a time to avoid getting webs with admissible 
non-canonical colorings.  

 There is a slight ambiguity in the rules above. At some point we
may have several choices which yield different cobordisms, depending on the
order in which we build them up. To remove this ambiguity we order all arcs of
the link, which induces an ordering of all unmarked edges in any of its
resolutions,
By convention we first build up the square or digon which contains the lowest
order edge.

\begin{figure}[h]
\includegraphics[height=0.7in]{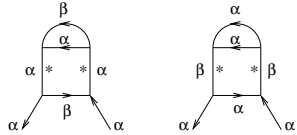}
\caption{Square-digons with coloring}
\label{fig:sqdigon_col}
\end{figure}

With these two observations in mind, it is not hard to see that the rules defining 
$\Sigma_{\phi}(\Gamma)$ are consistent. One can check that the two different ways of 
defining it for the square-digon webs in 
figure~\ref{fig:sqdigon_col}
yield the same cobordism indeed.

Recall that the definition of $R_{\phi}(\Gamma)$, which appears 
in the following lemma, can be found in the proof of 
lemma~\ref{lem:decomp2}, where we showed that $R_{\phi}(\Gamma)\cong 
Q_{\phi}(\Gamma)R(\Gamma)$. 

\begin{lem}
\label{lem:cyclic2} 
$C_{\phi}(\Gamma)$ is a free cyclic 
$R_{\phi}(\Gamma)$-module generated by $\Sigma_{\phi}(\Gamma)$, for any 
$\phi\in S(\Gamma)$. 
\end{lem}
\begin{proof} 
We use induction with respect to $v$, the number of trivalent vertices in 
$\Gamma$. 
The claim is obviously true for a circle. Suppose $\Gamma$ has a digon, 
with the edges ordered as in figure~\ref{fig:edge_order_digon}. Note 
that $X_1=X_4\in R(\Gamma)$ holds as a consequence of the relations 
in (\ref{eq:relations}). Let 
$\Gamma'$ be the web obtained by removing the digon, as in lemma~\ref{lem:KhK}.  
The possible canonical colorings of $\Gamma$ and the corresponding canonical colorings of $\Gamma'$ are 
$$\xymatrix@R=1.8mm{
\includegraphics[height=0.7in]{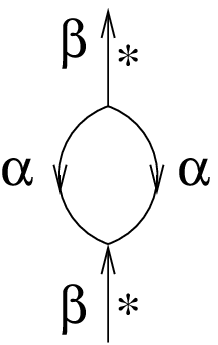} &
\includegraphics[height=0.7in]{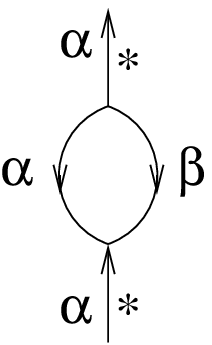} &
\includegraphics[height=0.7in]{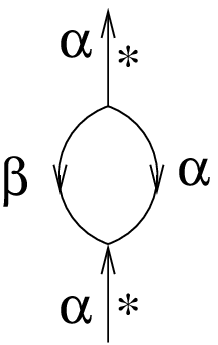} &
\includegraphics[height=0.7in]{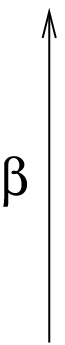} &
\includegraphics[height=0.7in]{edgeupa} &
\includegraphics[height=0.7in]{edgeupa} 
\\
\phi_1 & \phi_2 & \phi_3 & \phi'_1 & \phi'_2 & \phi'_3 
}.$$ 
We treat the case of $\phi_1$ first. 
Since the \emph{Digon Removal} isomorphism in lemma~\ref{lem:KhK} commutes with the action of $X_1=X_4$, 
we see that 
$$C_{\phi_1}(\Gamma)\cong C_{\phi_1'}(\Gamma')\oplus 
C_{\phi_1'}(\Gamma').$$
By induction $C_{\phi_1'}(\Gamma')$ is a free 
cyclic $R_{\phi_1'}(\Gamma')$-module generated by 
$\Sigma_{\phi_1'}(\Gamma')$. Note that the isomorphism 
maps 
$$\left(\Sigma_{\phi_1'}(\Gamma'),0\right)\quad\text{and}\quad 
\left(0,\Sigma_{\phi_1'}(\Gamma')\right)$$
to 
$$\Sigma_{\phi_1}(\Gamma)\quad\text{and}\quad 
X_2\Sigma_{\phi_1}(\Gamma).$$ Since 
$\dim R_{\phi_1}(\Gamma)=2\dim R_{\phi_1'}(\Gamma')$, we conclude that 
$C_{\phi_1}(\Gamma)$ is a free cyclic $R_{\phi_1}(\Gamma)$-module 
generated by $\Sigma_{\phi_1}(\Gamma)$.

Now, let us consider the case of $\phi_2$ and $\phi_3$. The \emph{Digon Removal} isomorphism in 
lemma~\ref{lem:KhK} yields
$$C_{\phi_2}(\Gamma)\oplus C_{\phi_3}(\Gamma)
\cong C_{\phi_2'}(\Gamma')\oplus 
C_{\phi_3'}(\Gamma').$$
Note that $\phi_2'=\phi_3'$ holds and by induction 
$C_{\phi_2'}(\Gamma')=C_{\phi_3'}(\Gamma')$ 
is a free cyclic $R_{\phi_2'}(\Gamma')=R_{\phi_3'}(\Gamma')$-module.  
As in the previous case, by definition of the canonical generators, it is easy 
to see that the isomorphism maps 
$$R_{\phi_2'}(\Gamma')\Sigma_{\phi_2'}(\Gamma')\oplus 
R_{\phi_3'}(\Gamma')\Sigma_{\phi_3'}(\Gamma')$$
to 
$$R_{\phi_2}(\Gamma)\Sigma_{\phi_2}(\Gamma)\oplus 
R_{\phi_3}(\Gamma)\Sigma_{\phi_3}(\Gamma).$$
Counting the dimensions on both sides of the isomorphism, we see that this proves 
the claim in the lemma for $C_{\phi_2}(\Gamma)$ and 
$C_{\phi_3}(\Gamma)$.

We could also have
$$\raisebox{-0.32in}{\includegraphics[height=0.7in]{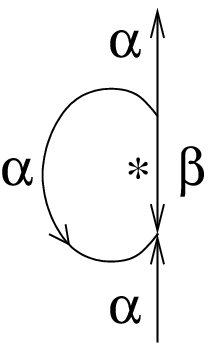}}\,,$$
but the same arguments as above apply to this case.

If we have several digons as in figure~\ref{fig:sev-dig}, similar arguments prove the 
claim when we stick in two digons at a time. 

Next, suppose $\Gamma$ contains a square, with the edges ordered as 
in figure~\ref{fig:edge_order_square} left. 
Let $\Gamma'$ and $\Gamma''$ be the two corresponding webs under the 
\emph{Square Removal} isomorphism 
in lemma~\ref{lem:KhK}. There is a number 
of possible canonical colorings. Note that there is no canonical coloring which colors 
all external edges by $\beta$. To prove the claim for all canonical 
colorings it suffices to consider only two: the one in which all external 
edges are colored by $\alpha$ and the coloring 
$$\includegraphics[height=0.7in]{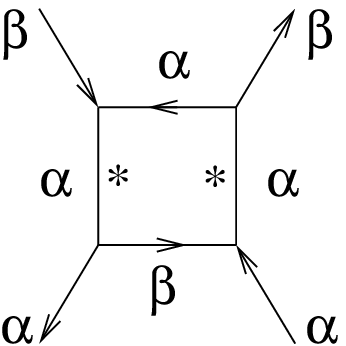}$$
All other cases are 
similar. 
Suppose that all external edges are colored by $\alpha$, then there are 
two admissible colorings:
$$\xymatrix@R=1.8mm{
\includegraphics[height=0.7in]{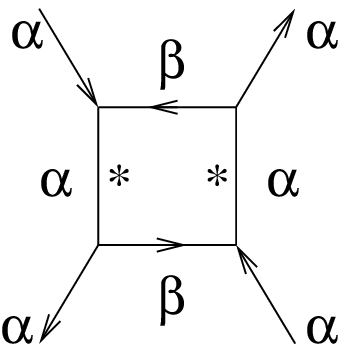} & \text{and}&
\includegraphics[height=0.7in]{sqcol_aaaababa} 
\\
\phi_1 && \phi_2
}.$$
Note that only $\phi_2$ can be canonical. 
Clearly there are unique canonical colorings of $\Gamma'$ and 
$\Gamma''$, with both edges colored by $\alpha$, 
which we denote $\phi'$ and $\phi''$. The isomorphism yields
$$C_{\phi_1}(\Gamma)\oplus C_{\phi_2}(\Gamma)\cong 
C_{\phi'}(\Gamma')\oplus C_{\phi''}(\Gamma'').$$
Suppose that the two edges in 
$\Gamma'$ belong to the same $\alpha$-cycle. 
We denote the number of $\alpha$-cycles 
in $\Gamma'$ by $m$. Note that the number of $\alpha$-cycles in 
$\Gamma''$ equals 
$m+1$. By induction $C_{\phi'}(\Gamma')=
R_{\phi'}(\Gamma')\Sigma_{\phi'}(\Gamma')$ and 
$C_{\phi''}(\Gamma'')=
R_{\phi''}(\Gamma'')\Sigma_{\phi''}(\Gamma'')$ 
are free cyclic modules of dimensions $2^m$ and 
$2^{m+1}$ respectively. Since $\phi_1$ is non-canonical, we know that 
$C_{\phi_1}(\Gamma)\cong C_{\phi'}(\Gamma')$, using the isomorphisms in 
(\ref{eq:canonical2}) and the results above about digon-webs. 
Therefore, we see that $\dim C_{\phi_1}(\Gamma)=2^m$ and 
$\dim C_{\phi_2}(\Gamma)=2^{m+1}$. By construction, we have 
$$\Sigma_{\phi''}(\Gamma'')\mapsto \left(*,
\Sigma_{\phi_2}(\Gamma)\right).$$
The isomorphism commutes with the actions on the external edges and 
$R_{\phi_2}(\Gamma)$ is isomorphic to $R_{\phi''}(\Gamma'')$, so we 
get 
$$R_{\phi_2}(\Gamma)\Sigma_{\phi_2}(\Gamma)\cong 
R_{\phi''}(\Gamma'')\Sigma_{\phi''}(\Gamma'').$$
For dimensional reasons, this implies that $C_{\phi_2}(\Gamma)$ is a 
free cyclic $R_{\phi_2}(\Gamma)$-module generated by $\Sigma_{\phi_2}(\Gamma)$. 

Now suppose that the two edges in $\Gamma'$ belong to different 
$\alpha$-cycles. This time we denote the number of $\alpha$-cycles in 
$\Gamma'$ and $\Gamma''$ by $2^{m+1}$ and $2^m$ respectively. 
We still have $C_{\phi_1}(\Gamma)\cong C_{\phi'}(\Gamma')$, so, by induction, 
we have $\dim C_{\phi_1}(\Gamma)=2^{m+1}$ and, therefore, $C_{\phi_2}(\Gamma)=2^m$. 
Consider the intermediate web $\Gamma'''$ colored by 
$\phi'''$ as indicated and the map between $\Gamma$ and $\Gamma'''$ in 
figure~\ref{fig:interm_square}. 
\begin{figure}[h]
$$\xymatrix@C=12mm@R=1mm{
\includegraphics[height=0.7in]{sqcol_aaaababa}  \ar[r]&
\includegraphics[height=0.7in]{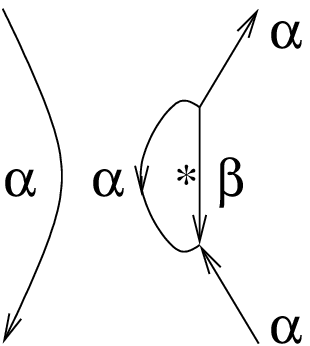}  \\
 \Gamma & \Gamma'''
}$$
\caption{}
\label{fig:interm_square}
\end{figure}
By induction, $C_{\phi'''}(\Gamma''')$ 
is a free cyclic $R_{\phi'''}(\Gamma''')$-module generated by 
$\Sigma_{\phi'''}(\Gamma''')$. By construction, 
we see that 
$\Sigma_{\phi_2}(\Gamma)$ is mapped to 
$$X_6\Sigma_{\phi'''}(\Gamma''')-
X_1\Sigma_{\phi'''}(\Gamma'''),$$ 
which is non-zero. Similarly we see that 
$X_1\Sigma_{\phi_2}(\Gamma)$ is mapped to 
$$
\begin{array}{l}
X_1X_6\Sigma_{\phi'''}(\Gamma''')-
X_1^2\Sigma_{\phi'''}(\Gamma''')=\\
X_1X_6\Sigma_{\phi'''}(\Gamma''')-
(2\alpha X_1-\alpha^2)\Sigma_{\phi'''}(\Gamma''').
\end{array}
$$ 
The latter is 
also non-zero and linearly independent from the first element. 
Since the map clearly commutes with the action of all elements 
not belonging to edges in the $\alpha$-cycle of $X_1$, the above shows that,  
for any nonzero element $Z\in R_{\phi_2}(\Gamma)$, the image of 
$Z\Sigma_{\phi_2}(\Gamma)$ in $C_{\phi'''}(\Gamma''')$ is non-zero.  
Therefore, we see that 
$$\dim R_{\phi_2}(\Gamma)\Sigma_{\phi_2}(\Gamma)=
\dim R_{\phi_2}(\Gamma)=2^m.$$
For dimensional reasons we conclude that 
$Q_{\phi_2}(\Gamma)C(\Gamma)$ is a free cyclic $R_{\phi_2}(\Gamma)$-module 
generated by $\Sigma_{\phi_2}(\Gamma)$. 

Finally, let us consider the canonical coloring 
$$\includegraphics[height=0.7in]{sqcol_baababaa}\raisebox{22pt}{.}$$
In this 
case $C_{\phi}(\Gamma)$ 
is isomorphic to $C_{\phi''}(\Gamma'')$, where $\phi''$ is the 
unique compatible canonical coloring of $\Gamma''$, because there is no 
compatible coloring of $\Gamma'$. Note that $R_{\phi}(\Gamma)$ is isomorphic 
to $R_{\phi''}(\Gamma'')$ and $\Sigma_{\phi}(\Gamma)$ is 
mapped to $\Sigma_{\phi''}(\Gamma'')$. Therefore the result 
follows by induction.
\end{proof}

Finally we arrive at our main theorem in this subsection. 

\begin{thm}
\label{thm:finalformula2}
$$\uht i L\cong \bigoplus_{L'\subseteq L} \kh {i-j(L')} * {L'},$$
where $j(L')=2\lk(L',L\backslash L')$.
\end{thm}
\begin{proof}

By theorem~\ref{thm:decomp2} we know that 
$$\uht i L=\bigoplus_{\phi\in S(L)}\uht i {L_{\phi}}.$$
Let $\phi\in S(L)$ be fixed and let $L_{\alpha}$ be the sublink of $L$ 
consisting of the components colored by $\alpha$. We claim that 
\begin{equation}
\label{eq:finalformula2}
\uht i {L_{\phi}}\cong \kh {i-j(L')} * {L_{\alpha}},
\end{equation}
from which the theorem follows. First note that, without loss of generality, 
we may assume that $\alpha=0$, because we can always apply the isomorphism 
$\planedot\mapsto\planedot - \alpha\plane$. 
Let $C_{\phi}(\Gamma)$ be a canonical cochain group.  
By lemma~\ref{lem:cyclic2}, we know that $C_{\phi}(\Gamma)$ is a free 
cyclic $R_{\phi}(\Gamma)$-module generated by 
$\Sigma_{\phi}(\Gamma)$. Therefore we can identify any 
$X\in R_{\phi}(\Gamma)$ with $X\Sigma_{\phi}(\Gamma)$. 
There exist isomorphisms 
\begin{equation}
\label{eq:iso}
R_{\phi}(\Gamma)\cong\bC[X_i\,|\,\phi(i)=\alpha]/\left(
X_i+X_j, X_i^2\right)\cong A^{\otimes m},
\end{equation}
where $A=\bC[X]/\left(X^2\right)$.
As before, the relations $X_i+X_j=0$ hold whenever the edges $i$ and $j$ 
share a common trivalent vertex and $m$ is the number of $\alpha$-cycles in 
$\Gamma$. Note also that $X_iX_j=0$ holds, if $i$ and $j$ share a trivalent 
vertex. The first isomorphism in (\ref{eq:iso}) is immediate, but for 
the definition of the second isomorphism we have to 
make some choices. First of all we have to chose an ordering of the arcs of 
$L$. This 
ordering induces a unique ordering on all the unmarked edges of $\Gamma$, where 
we use Bar-Natan's~\cite{bar-natan} convention that, if an edge in $\Gamma$ is 
the fusion of two arcs of $L$, we assign to that edge the smallest of the two 
numbers. Now delete all edges colored by $\beta$. 
Consider a fixed $\alpha$-cycle. In this $\alpha$-cycle 
pick the edge $i$ with the smallest number 
in our ordering. This edge has an orientation induced by the orientation of 
$L$.

We identify the $\alpha$-cycle with a circle, by deleting all vertices in 
the $\alpha$-cycle, oriented according to the 
orientation of the edge $i$. If the circle is oriented clock-wise we say that 
it is negatively oriented, otherwise we say that it is positively oriented. 
The circles corresponding to the $\alpha$-cycles are ordered 
according to the order of their minimal edges. They can be 
nested. As in Lee's paper~\cite{lee} we say that a circle is positively nested 
if any ray from that circle to infinity crosses the other circles 
in an even number of points, otherwise we say that it is negatively nested.    
The isomorphism in (\ref{eq:iso}) is now defined as follows. Given the 
$r$-th $\alpha$-cycle with minimal edge $i$ we define  
$$X_i\mapsto \epsilon 1\otimes \cdots\otimes X\otimes\cdots
\otimes1,$$
where $X$ appears as the $r$-th tensor factor. If the 
orientation and the nesting of the $\alpha$-cycle have the same sign, then 
$\epsilon=+1$, and if the signs are opposite, then $\epsilon=-1$.  
The final result, i.e. the claim of this theorem, holds true no matter 
which ordering of the arcs of $L$ we begin with.  
\begin{figure}[h]
$$
\includegraphics[height=0.2in]{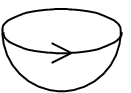}\; 
\includegraphics[height=0.2in]{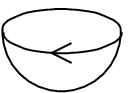}\;\raisebox{8pt}{$\mapsto$}
\;\raisebox{-0.1in}{\includegraphics[height=0.3in]{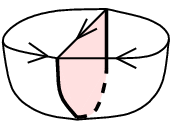}}\quad
\raisebox{8pt}{,}\quad
\raisebox{-0.1in}{\includegraphics[height=0.3in]{halftheta}}
\;\raisebox{8pt}{$\mapsto$}\;
\includegraphics[height=0.2in]{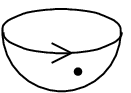}\; 
\includegraphics[height=0.2in]{cuplor}\;
\raisebox{8pt}{-}\;
\includegraphics[height=0.2in]{cupror}\; 
\includegraphics[height=0.2in]{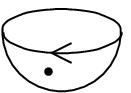}\;
$$

$$\begin{matrix}
\includegraphics[height=0.35in]{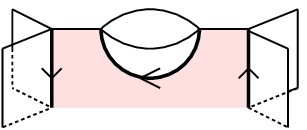} &\raisebox{10pt}{$\mapsto$}&
\includegraphics[height=0.35in]{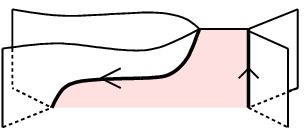}
\\
\includegraphics[height=0.35in]{dedge_id} &\raisebox{12pt}{$\mapsto$}&
\includegraphics[height=0.35in]{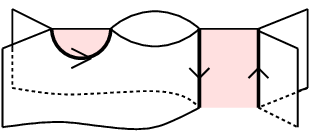} & \raisebox{12pt}{=} &
\includegraphics[height=0.35in]{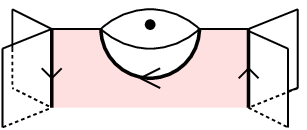} &\raisebox{12pt}{-}&
\includegraphics[height=0.35in]{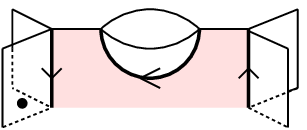}
\end{matrix}$$

$$
\includegraphics[height=0.6in]{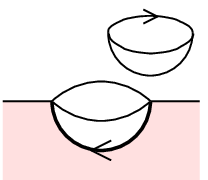}\;\raisebox{10pt}{$\mapsto$}
\includegraphics[height=0.6in]{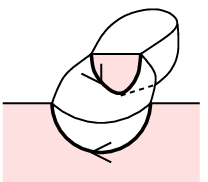}\quad
\raisebox{8pt}{,}\quad
\includegraphics[height=0.6in]{scupscupout}\;\raisebox{10pt}{$\mapsto$}
\includegraphics[height=0.6in]{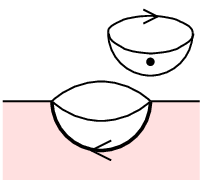}\;\raisebox{10pt}{-}\;\includegraphics[height=0.6in]{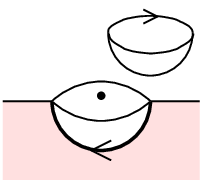}
$$

$$
\left.
\begin{matrix}
\raisebox{-0.3in}{\includegraphics[height=0.7in]{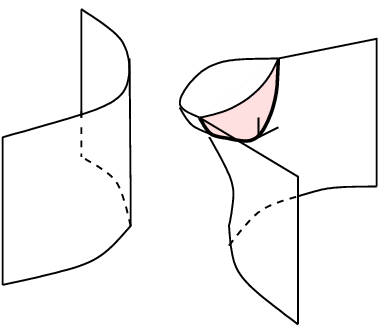}}\mapsto  
\raisebox{-0.3in}{\includegraphics[height=0.7in]{ssquareupr}} 
   \\  \\
\raisebox{-0.3in}{\includegraphics[height=0.7in]{ssquareupr}}\mapsto 
\raisebox{-0.3in}{\includegraphics[height=0.7in]{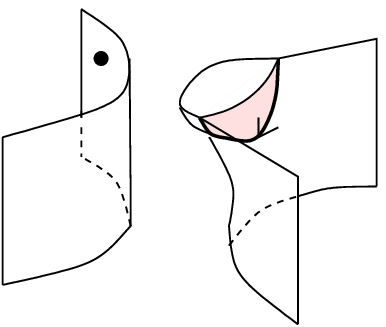}}
 -  
\raisebox{-0.3in}{\includegraphics[height=0.7in]{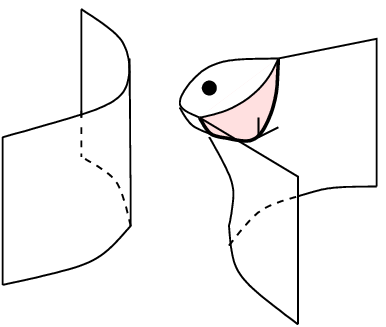}}
\end{matrix}
\right\}
\begin{matrix}
\raisebox{-0.2pt}{\includegraphics[height=3.2mm]{tedgesop}}\;\text{ belong to same}
\\ \\
\text{$\alpha$-cycle in }\raisebox{-0.3in}{\includegraphics[height=0.6in]{sqcol_aaaababa}}
\end{matrix}
$$

$$
\left.
\begin{matrix}
\raisebox{-0.3in}{\includegraphics[height=0.7in]{ssquareup}}\mapsto 
\raisebox{-0.3in}{\includegraphics[height=0.7in]{sedgedigonr}}
   \\ \\
\raisebox{-0.3in}{\includegraphics[height=0.7in]{sedgedigonr}}\mapsto  
\raisebox{-0.3in}{\includegraphics[height=0.7in]{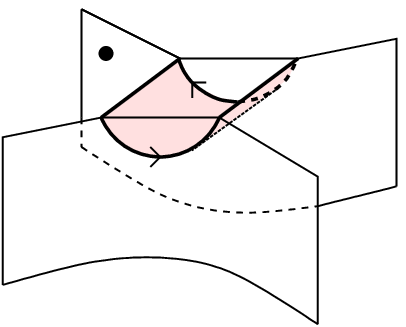}}+
\raisebox{-0.3in}{\includegraphics[height=0.7in]{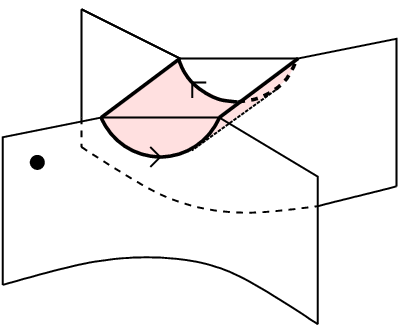}}
\end{matrix}
\right\}
\begin{matrix}
\raisebox{-0.2pt}{\includegraphics[height=3.2mm]{tedgesop}}\;\text{ belong to different} 
\\ \\
\text{$\alpha$-cycles in }\raisebox{-0.3in}{\includegraphics[height=0.6in]{sqcol_aaaababa}}
\end{matrix}
$$

\caption{Behaviour of canonical generators under elementary cobordisms}
\label{fig:cgen_beh_elemcob}
\vspace*{0.6mm}
\end{figure}
It is easy to work out the behaviour of the canonical generators with respect 
to the elementary cobordisms as can be seen in 
figure~\ref{fig:cgen_beh_elemcob}. 
For the cobordisms shown in figure~\ref{fig:cgen_beh_elemcob} having one or more cycles there is also a version with one cycle inside the other cycle or 
a cycle inside a digon.  

The two bottom maps in figure~\ref{fig:cgen_beh_elemcob} require some 
explanation. Both can only be understood by considering all possible closures 
of the bottoms and sides of their sources and targets. Since, by definition, 
the canonical generators are constructed step by step introducing the 
vertices of the webs in some order, we can assume, without loss of generality, 
that the first vertices in this construction are the ones shown. With 
this assumption the open webs at the top and bottom of the cobordisms in 
the figure are to be closed only by simple curves, without vertices, and 
the closures of these cobordisms, outside the bits which are shown, only use 
cups and identity cobordisms. Baring this in mind, the claim implicit in 
the first map is a consequence of relation~\ref{eq:usefuleq}. 

For the second map, recall that \includegraphics[height=2.3mm]{tedgesop} 
belong to different $\alpha$-cycles. Therefore there are two different ways 
to close the webs in the target and source: two cycles side-by-side or one  
cycle inside another cycle. We notice that from theorem~\ref{thm:decomp2} we 
have the isomorphism
$$\raisebox{-0.25in}{\includegraphics[height=0.6in]{sqcol_aaaababa}}\cong
\raisebox{-0.25in}{\includegraphics[height=0.6in]{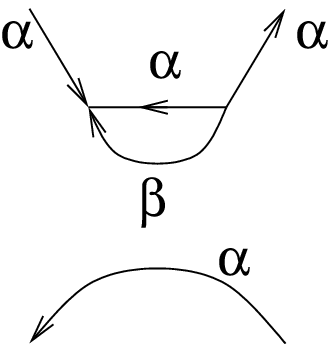}}.$$
We apply this isomorphism to the composite of the source foam and 
the elementary foam and to the target foam of the last map in 
figure~\ref{fig:cgen_beh_elemcob}. Finally use equation~\ref{eq:usefuleq} and 
relation (\emph{CN}) to the former to see that both foams are isotopic.

Note that in $\cut * {L_{\phi}}$ we only have to consider elementary 
cobordisms at crossings between two strands which are both colored by 
$\alpha$. With the identification of $R_{\phi}(\Gamma)$ and 
$A^{\otimes m}$ as above, it is now easy to see that 
the differentials in $\cut * {L_{\phi}}$ behave exactly as in Khovanov's 
original $sl_2$-theory for $L_{\alpha}$.  

The degree of the 
isomorphism in (\ref{eq:finalformula2}) is easily computed using the fact 
that in both theories the oriented resolution has homological degree zero. 
Therefore we get an 
isomorphism
$$\uht i {L_{\phi}}\cong \kh {i-j(L')} * {L_{\alpha}}.$$   
\end{proof}


\vspace*{1cm}

\noindent {\bf Acknowledgements} We thank Mikhail Khovanov and Sergei 
Gukov for enlightening conversations and exchanges of email.

The first author was supported by the 
Funda\c {c}\~{a}o para a Ci\^{e}ncia e a Tecnologia through the
programme ``Programa Operacional Ci\^{e}ncia, Tecnologia, Inova\c
{c}\~{a}o'' (POCTI), cofinanced by the European Community fund FEDER.


\end{document}